\pdfoutput=1
\documentclass{article}
\usepackage{amsfonts}
\usepackage{amsmath}
\usepackage{mathtools}
\usepackage{graphicx}
\usepackage[normalem]{ulem}
\usepackage{url}
\usepackage[utf8]{inputenc}
\usepackage{bbm}
\usepackage{algorithm}
\usepackage{caption}
\usepackage{subcaption}
\usepackage{algpseudocode}
\usepackage{multicol}
\usepackage{enumitem} 
\usepackage{wasysym} 
\usepackage{hyperref}
\usepackage[dvipsnames,svgnames]{xcolor}
\usepackage{authblk}
\newcommand{\indic}{\mathbbm{1}}
\newcommand{\mean}{\mathbb{E}}

\newcommand{\prob}{\mathbb{P}}

\graphicspath{{figures/}}

\title{Coupling and selecting constraints in Bayesian optimization under uncertainties}
\date{}

\author[1]{Julien Pelamatti}
\author[2]{Rodolphe Le Riche}
\author[3]{Céline Helbert}
\author[3]{Christophette Blanchet-Scalliet}
\affil[1]{\small EDF R\&D, France (formerly Ecole des Mines de Saint-Etienne)}
\affil[2]{CNRS LIMOS (Mines Saint-Etienne and UCA), France}
\affil[3]{Université de Lyon - CNRS, UMR 2508, Institut Camille Jordan - Ecole Centrale de Lyon, France}
\begin{document}

\normalsize 
\maketitle

\begin{abstract}
\normalsize 
We consider chance constrained optimization where it is sought to optimize a function while complying with constraints, both of which are affected by uncertainties. 
The high computational cost of realistic simulations strongly limits the number of evaluations and makes this type of problems particularly challenging.
In such a context, it is common to rely on Bayesian optimization algorithms. 
Assuming, without loss of generality, that the uncertainty comes from some of the inputs, 
it becomes possible to build a Gaussian process model in the joint space of design and uncertain variables. 
A two-step acquisition function is then used to provide, both, promising optimization variables associated to relevant uncertain samples.
Our overall contribution is to correlate the constraints in the GP model and exploit this to optimally decide, at each iteration, which constraint should be evaluated and at which point.
The coupled Gaussian model of the constraints relies on an output-as-input encoding. The constraint selection idea is developed by enabling that each constraint can be evaluated for a different uncertain input, thus improving the refinement efficiency. Constraints coupling and selection are gradually implement in 3 algorithm variants which are compared to a reference Bayesian approach.
The results are promising in terms of convergence speed, accuracy and stability as observed on a 2, a 4 and a 27-dimensional problems.

\end{abstract}

\section{Introduction}

Models of realistic physical or engineering systems are frequently impacted by uncertain random parameters that represent model error, uncontrolled environmental conditions, or numerical noise such as imperfectly converged internal routines. 
When optimizing these systems, it is necessary to take into account the presence of the random parameters to determine a sufficiently robust solution in terms of both the objective and the constraint functions. 
Depending on the application, different approaches to handle the uncertainties within an optimization problem exist. 

A distinction can be made between the stochastic and the robust approaches. 
In the stochastic methods, the uncertainty in the optimization criteria (objective
and constraint functions) is modeled as a random noise corrupting the observed values. 
Stochastic algorithms can, under conditions about their own stochasticity, accomodate the noise in the observations and converge to the problem solutions: 
this has given rise in the early 50's to stochastic descent algorithms \cite{kiefer1952stochastic,arrow1958studies} that have since then experienced great developments \cite{spall2005introduction,andrieu2011gradient} often in relation to machine learning \cite{kingma2014adam}. 
Well-performing versions of (stochastic) evolutionary algorithms for noisy functions have been identified \cite{brockhoff_miror_2010,loshchilovsaACM2012} thanks to competitions on reference benchmarks \cite{COCO2012noisy}.

Other approaches to optimization under uncertainty have emerged in many communities (operations research \cite{gabrel2014recent}, finance \cite{fabozzi2007robust}, engineering \cite{youn2016structural}, machine learning \cite{caramanis201214} and many others). They quantify with deterministic functions the trade-off between the performance and the risk associated to the uncertainties. In robust
optimization, the set of uncertain parameters is typically handled through
worst-case scenarii or convex approximations of it  \cite{gabrel2014recent,ben2009robust}. 
In a variety of situations where
the set of uncertainties is large (possibly infinite) and the problem structure non-convex, over-conservatism of the solutions and the computational tractability
of the methods are inherent difficulties of this family of problems. 
When the uncertainties are described by probability distributions, it is now customary to optimize a risk measure which takes into account the whole probability distribution. 

In the presence of optimization constraints, the uncertainties affect the feasibility of the solutions. 
Probabilistic models of the constraints are called chance constraints \cite{nemirovski2012safe} or reliability constraints \cite{valdebenito2010survey}.
The handling of chance constraints, mostly in engineering, has spawned the field of Reliability-Based Design Optimization (RBDO). 
In RBDO, a popular family of approaches appproximate the probability of constraints satisfaction by reliability indices which are computationally easier to calculate and suitable to a coupling within an optimization problem (First/Second Order Reliability Methods, \cite{valdebenito2010survey}). 
However, when working with reliability indices, the failure probability is often poorly estimated. 
Other RBDO techniques approximate the failure probability with subset and importance sampling in conjunction with metamodels which allow the inclusion of the chance constraint within an optimization loop \cite{dubourg2011reliability,moustapha2016quantile}.

In this article, the problem of minimizing the mean of a stochastic function (also known as Bayes risk) under chance constraints (i.e., reliability constraints) is addressed. 
It is assumed that the problem functions are costly (i.e., cannot be calculated more than $\mathcal O(100)$ times) so that one seeks to find the neighborhood of the optimal solution with as few function evaluations as possible. 
Examples of expensive functions are nonlinear Finite Element Models used for structural analyses, Computational Fluid-Dynamics simulations and most real life experiments. 
The cost limitation rules out the possibility of relying on most stochastic algorithms (like stochastic gradients, evolution strategies) which require a large number of evaluations. 
It is also assumed that there is no information regarding the convexity or the gradient of the problem functions. 
These assumptions limit the type of algorithms which are relevant to the class of surrogate-based methods, among which Bayesian Optimization (BO) \cite{jones1998efficient,mockus2012bayesian} will be our focus.

The principle underlying BO is to learn the problem functions by Gaussian processes (GPs) \cite{rasmussen2003gaussian} which are computationally less expensive than the functions. An acquisition criterion, based on the GPs and adapted to the problem formulation, defines the most promising points of the design space. Its repeated optimization, followed by an update of the GPs, constitutes the BO algorithm. 

\vskip\baselineskip
Only recently has the handling of chance constraints been treated with BO \cite{dubourg2011reliability,moustapha2016quantile,beland2017bayesian,cousin2020} (cf. Section~\ref{sec:pastRelatedWorks} for a detailed discussion of past related works).
This article presents three contributions to BO in the situation where several inequalities must be satisfied in probability.
Similarly to \cite{beland2017bayesian,elamri2021sampling}, the starting point is a procedure which first maximizes the feasible improvement and then reduces its variance in order to define new values of both the optimization and the random variables.

The first contribution of this article is to account for the couplings between the constraints through an ``output as input'' formulation. 
This is important as the constraints are usually linked to each other, e.g., when they express different criteria extracted from the same physical phenomenon, or when they are competing anticorrelated criteria (such as price vs. availability of raw materials).

The second contribution concerns the common situation where the objective and the constraints are calculated with different codes. 
We will show how the objective and the constraints can be refined independently in the random search space.
Indeed, the random variables that provide the largest amount of information regarding \textit{i)} the expected objective function and \textit{ii)} the expected feasibility of a given candidate solution are different.

The third contribution is to allow the optimization algorithm to select, at each iteration, the constraint that provides the most information about the location of the robust optimum. The motivation is that, when the constraints can be computed independently, it is useless to refine a constraint for which the (in)feasibility of the candidate solution is known with a sufficiently high certainty.

The article continues with a more detailed review of past related works followed, in Section~\ref{ProbForm}, by a formulation of the problem. 	
Section~\ref{Multioutput} describes the output as input model and the BO acquisition criteria. 
The independent calculation of the random variables and the selection of the constraint to be refined are presented in Section~\ref{ProxExt}.
The last parts of the article are dedictated to the implementation of the algorithms (Section~\ref{ImplemNotes}) and numerical tests (Section~\ref{NumRes}) on 2 analytical and an industrial problems. 
To assist the reader, the main symbols and acronyms used in the article are summarized in Appendix~\ref{sec:notations}.

\section{Past related works}
\label{sec:pastRelatedWorks}

Our study rests on preliminary works concerning multi-output Gaussian processes, Bayesian optimization with constraints and Bayesian optimization with uncertainties. 
We now review this bibliography in turn.

\subsection{Multi-output Gaussian Processes}
A contribution of this article is the modeling of chance constraints with a single multi-output Gaussian process rather than with independent models. The idea is to take into account the correlation between outputs in order to improve the modeling accuracy thus reducing the number of computer runs required to converge to the problem optimum.

The model we introduce relies on the ``output as input'' principle where a nominal input is added to the variables list to designate the selected constraint. 
Other  approaches found in the litterature can also model correlations between outputs. A first parallel can be drawn with surrogates for multi-fidelity modeling \cite{peherstorfer2018survey}, \cite{fernandez2016review}. 
These methods represent with a single vector valued function the responses of a physical system associated to different precisions (\emph{e.g.}, varying mesh densities in finite elements). 

Most of the so-called co-kriging approaches (\cite{kennedy2000predicting}, \cite{perdikaris2017nonlinear}, \cite{le2014recursive}) define a hierarchical relation between the various fidelity levels:
the surrogate model for the fidelity level $t$ is recursively defined as the Gaussian process model of the previous fidelity level ($t-1$) summed with a second Gaussian process, modeling the additional information provided by the next fidelity level.
However, such a hierarchical dependence hypothesis doesn't always stand. 
For example, when the components of the vectorial output represent different physical quantities, it doesn't necessarily make sense to model them recursively. 
In this case, the order of the various outputs would be an arbitrary modeling decision.

Alternative, non hierarchical GPs for vector valued functions include models based on co-regionalization matrices between the outputs \cite{williams2007multi}, \cite{goulard1992linear}, \cite{alvarez2011kernels}. 
Recently, deep Gaussian processes have also been used to model multi-fidelity phenomena, allowing to better capture non-stationary behaviors (\cite{brevault2020overview}, \cite{cutajar2019deep}).

The multi-output modeling methods described above have been applied to computationally intensive optimization problems. 
These include multi-objective optimization (\cite{shah2016pareto}, \cite{ghoreishi2019multi}) and multi-fidelity design optimization (\cite{meliani2019multi}, \cite{brooks2011multi}, \cite{garland2020multi}).
The ability to select which constraint should be evaluated at each iteration in order to maximise the information gained is a salient feature of the current 
article. 
A somewhat comparable idea was proposed in \cite{tao2021multi} to deal with multi-model design. A specific sub-model (and its intervening variables) was selected for evaluation based on a ratio of the ``knowledge gradient'' (the predicted one-step-ahead gain in performance) to the submodel cost.

\subsection{Constrained Bayesian Optimization}
The current article handles the chance constraints by aggregating them with the objective function through the expected feasible improvement acquisition criterion that was first proposed in \cite{schonlau1998global}.
Other techniques to account for constraints in Bayesian optimization exist. 

In \cite{sasena2002exploration}, the constraints are metamodeled but they are not aggregated with the objective function. 
The acquisition of a new point is a constrained optimization problem involving only the metamodels of the objective function and the constraints. This problem, which is not as costly as the original one, is treated with a nonlinear constrained optimization algorithm. 

Augmented Lagrangians have been adapted to Bayesian optimization in the presence of equality and inequality constraints in \cite{picheny2016bayesian}. Slack variables are introduced to transform the inequalities into equalities. An expected improvement formulated in terms of augmented Lagrangian values is given.

In \cite{audet2000surrogate}, a constrained Bayesian optimization is carried out inside a filter-based generalized pattern search algorithm.
Candidate points are created by the sequential application of two acquisition criteria, an expected violation of the constraints that must be below a threshold followed by an expected improvement.

\subsection{Bayesian optimization under uncertainty}
It is no wonder that Bayesian optimization has often been adapted to account for uncertainties. When the uncertainties are given as probability densities and the performance is quantified by risk measures, Gaussian processes, which are probabilistic mathematical objects, come as a natural approach.
The adaptation of BO to uncertainties affects either the Gaussian process (the parameters over which it is defined, the function that it models) and/or the acquisition criterion. 

\paragraph{Robust optimization.}
Many previous works have focused on the robust optimization of a function in the absence of constraints. 
When dealing with noisy objective functions, nugget effects have been dimensionned in proportion to the noise and added to the Gaussian process \cite{leriche2009gears}. 
The expected improvement acquisition criterion has been modified in various ways to account for this noise \cite{picheny2013benchmark}. 
New acquisition criteria have been proposed to accomodate this noise. 
Risk adverse criteria are based on quantiles \cite{picheny2013quantile,cakmak2020bayesian,sabater2021bayesian} while a risk neutral criterion is the correlated knowledge gradient \cite{scott2011correlated}. 
A risk adverse approach is described in \cite{torossian2020bayesian} that encompasses both quantiles and expectiles:
specific Gaussian processes are given for the quantiles and expectiles along with two BO strategies, a generalization of the upper-confidence bound and a batch sampling strategy. The compromises between risk adverse and risk neutral criteria have been explored through a multi-objective approach in \cite{RIBAUD2020106913}.

A key idea in the treatment of uncertainties is the \textit{joined space}, 
which is the tensorization of the deterministic space of the optimization variables with the space of the random parameters. 
It opens the possibility to define the Gaussian processes and the acquisition in the joined space.
At each iteration then, a BO can provide not only a new value for the optimization variables, but also a new sample for the random parameters. 
Such a procedure is a simultaneous optimization and optimal sampling, where the optimality of the sampling is in the sense of the acquisition criterion.
An early example is given in \cite{janusevskis_jogo_2012} where the mean objective function is minimized and the acquisition criterion is the one-step-ahead variance of the objective function at desirable locations.
The same mean objective function is considered, still in a joined space, in \cite{zhang2017sequential}. The variables and random parameters are samples of an adapted importance sampling density which is calculated with the GP and which accounts sequentially for the contribution to the objective function and to the uncertainty in the objective function.\\
In \cite{cakmak2020bayesian}, the acquisition criterion is a knowledge gradient applied to a quantile (a conditional-value-at-risk) of the objective function. Because of the quantile, this work bears a resemblance to \cite{picheny2013quantile,torossian2020bayesian} but in the joined space, therefore also providing random parameters values.

\paragraph{Optimization with reliability constraints. }
The estimation of reliability constraints, outside of any optimization, with the help of Gaussian processes has received a lot of attention \cite{menz2020adaptive,dubourg2014meta,bect2012sequential,mattrand2021adaptive}.
Such reliability analyses based on metamodels have then been included within optimization searches.
A recent survey and a global framework for RBDO can be found in \cite{moustapha2019surrogate}.  
Several of these contributions involve kriging metamodels built over the space of design variables \cite{lacaze2013reliability,li2016local,chen2015important}. 

More recently, like in the case of noisy objective functions, optimization with reliability constraints has been tackled through Gaussian processes built in a space that allows to select both the variables and the noise parameters. 
In \cite{dubourg2011reliability}, \cite{moustapha2016quantile} and \cite{cousin2020} because the noise directly affects some of the variables, a formulation in the 
\emph{augmented space} (not to be mistaken with the joined space) is developped: 
the new optimization variables are parameters of the noise density functions and the Gaussian processes are built in the space of the original design vriables. 
The acquisition criteria, which command the computationally intensive simulations, primarily aim at refining the estimation of the reliability constraints. 
In \cite{dubourg2011reliability}, the acquisition criterion is the probability of a point to be in the margin of uncertainty of the constraint. In \cite{moustapha2016quantile}, it is the deviation number (a ratio of the average and standard deviation of the kriging model of the constraint).
In both works, the distribution variables are updated independently of calls to the simulation by a direct iteration on the optimization problem.

\cite{beland2017bayesian} provides a complete Bayesian problem formulation close to ours: the mean objective function is minimized with reliability constraints; the Gaussian processes are defined over the joined space; the acquisition criterion stems from a similar two steps logic, the first step being a maximization of the feasible improvement to define the optimization variable and the second step being a maximization of the posterior (integrated) variance to define the noise parameters. However the constraints are independent.

\vskip\baselineskip
To the authors knowledge, this article is the first fully Bayesian formulation of the problem 
of minimizing the expected value of a function under chance constraints where the constraints are 
not independent, which allows to select which constraint needs to be evaluated.



\section{Problem formulation}
 \label{ProbForm}
We work in the joined space of the optimization variables, $\mathbf{x} \in \mathcal{S_X} \in \mathbb{R}^d$, over which the user has control (e.g., targeted dimensions, masses, pressures), and the random parameters $\mathbf{U}  \in \mathcal{S_U} \in \mathbb{R}^m$ (e.g., atmospheric parameters, manufacturing errors, numerical noise) for which a probability density function $\rho_U$ is known. 
Although the $\mathbf{U}$'s are not controlled in the real world, they can be chosen during the numerical simulations through their samples $\mathbf u$.
We assume indeed that the objective and the constraint functions depend solely on $\mathbf{x}$ and $\mathbf{u}$ and can be calculated with simulation codes that are deterministic.
The postulate of $\rho_U$ is a strong hypothesis that considerably simplifies the calculations. Less restrictive approaches to the uncertainty description would still require expert knowledge in the form of a parameterized family of probability density functions that must be selected. 

We consider the following optimization problem where $f(\cdot)$ is the objective function and where there are $l$ constraints $\{ g_1(\cdot),\dots, g_{l}(\cdot)\}$: 
\begin{equation}
\begin{split}
& \text{find } \quad \mathbf{x}^\star = \arg\min_{\mathbf{x} \in \mathcal{S_X} \subset \mathbb{R}^d} \mean_{\mathbf{U}}[f(x,U)] \\
& \text{such that } \quad \prob_{\mathbf{U}}(g_p(\mathbf x^\star,U) \leq 0, p = 1, \dots, l) \geq 1 - \alpha 
\end{split}
\label{PBDef}
\end{equation}
The objective is to minimize the mean of $f(\cdot)$ with respect to $\mathbf{U}$ (also called Bayes risk) while simultaneously satisfying the constraints with a probability above a given threshold $\alpha$ (called reliability or chance constraint).  
Other problem formulations that include uncertainties exist \cite{Balesdent2020}. 
Of course, the formulation conditions the characteristics of the optimal design. 
In our work, the emphasis is placed on constraint satisfaction with a guaranteed reliability. 
The effect of the uncertain variables on the objective function is taken in expectation because of its analytical and numerical tractability. 
The average performance is the appropriate formulation in the frequent situations where the objective $f$ will be measured many times. 
Examples include the life span of electric bulbs, the energy consumption of a machine, the speed of convergence of an algorithm.

The objective function is considered to be independent from the constraints, while the constraints are taken as mutually dependent. 
Other formulations with quantiles conditioned by feasibility have been described in \cite{pujol2009incertitude}, yet they remain an open challenge for costly problems. 
Problem \eqref{PBDef} can be rewritten by expressing the probability of feasibility as the expectation of a Bernoulli variable representing whether or not all of the constraints are satisfied:
\begin{eqnarray*}
\prob_{\mathbf{U}}(g_p(\mathbf{x},\mathbf{U}) \leq 0,~ p=1,\dots,l ) \geq 1 - \alpha & \Leftrightarrow & \mean_{\mathbf{U}}[ \indic_{\{g_p(\mathbf{x},\mathbf{U}) \leq 0,~ p=1,\dots,l \}}]  \geq 1 - \alpha,\\
 &\Leftrightarrow & 1 - \alpha -\mean_{\mathbf{U}}[\indic_{\{g_p(\mathbf{x},\mathbf{U}) \leq 0,~ p=1,\dots,l \}}]  \leq 0, \\
 &\Leftrightarrow & c(\mathbf{x}) \leq 0.
\end{eqnarray*}
The optimization problem becomes
\begin{equation}
\mathbf{x}^\star = \arg \min\limits_{\mathbf{x} \in \mathcal{S_X}} ~z(\mathbf{x}) ~\text{s.t.}~ c(\mathbf{x}) \leq 0~,
\label{newproblem}
\end{equation}
where $z(.) = \mean_{\mathbf{U}}[f(.,\mathbf{U})]$.

\section{Robust Bayesian optimization with a multi-output model of the constraints}
\label{Multioutput}
Problems~\eqref{PBDef} or \eqref{newproblem} have been treated with a robust Bayesian approach in \cite{elamri2021sampling} where it was 
assumed that the constraints were uncorrelated.
Such an assumption has two main consequences: first of all, the constraints can be modeled by independent GPs which are trained with separate data sets.
Second, the probability that all of the constraints are satisfied at a given location in the design space can be computed as the product of the probabilities of feasibility of each constraint. 

However, the independence assumption may be overly simplistic and may potentially limit the algorithm performance. 
Real constraints are often positively or negatively correlated to each other because they are computed from the same underlying physical phenomenon. An example is the correlation between the maximum displacement and the Von Mises stress on a structure in bending. 
Another example is the negative correlation between manufacturing cost and the quality (e.g., the dimensional accuracy) of a product. 
Improving the model accuracy in the context of the optimisation of the return of a portofolio under a constraint on its variance is also a relevant example. The expected return is maximized under a constraint on the return variance, but the variance of the portfolio depends on the correlation between its financial assets.

For the abovementioned reasons, in this article, we consider correlated constraints. 
Not only will our model be more realistic, but it will have the ability to share information between the constraints, thus improving the accuracy of the model. 

\subsection{Gaussian regression for a multi-output model of the constraints}
\label{sec:gaussianRegressionMulti}
Let's assume that the objective function $f$ and the constraints  $\mathbf{g}(\cdot) = \{ g_1(\cdot),\dots, g_{l}(\cdot)\} $ are realizations of a Gaussian process $F$ and a multi-output GP $\mathbf{G}$, respectively. 
$\mathbf{G}$ is furthermore assumed to be independent of $F$, and both $F$ and $\mathbf{G}$ are built in the joint space $\mathcal{S_X}\times \mathcal{S_U}$. As GPs, $F$ and $\mathbf{G}$ are characterized by a prior mean and a covariance function:
\begin{eqnarray}
F(\mathbf{x}, \mathbf{u})  &\sim & \mathcal{GP}(m_F(\mathbf{x}, \mathbf{u}) , k_F(\mathbf{x}, \mathbf{u},\mathbf{x}', \mathbf{u}')) \\
\mathbf{G}(\mathbf{x}, \mathbf{u}) = \left(
\begin{array}{c }
G_1(\mathbf{x}, \mathbf{u}) \\
\vdots  \\
G_l(\mathbf{x}, \mathbf{u}) 
\end{array} \right) &\sim & \mathcal{GP}(\mathbf{m}_\mathbf{G}(\mathbf{x}, \mathbf{u}), \mathbf{k}_ \mathbf{G}(\mathbf{x}, \mathbf{u},\mathbf{x}', \mathbf{u}'))
\end{eqnarray}
where $m_F$ and $\mathbf{m_G}$ are the mean functions and $k_F$ and $\mathbf{k_G}$ are the covariance functions. 
While $k_F$ is scalar valued, $\mathbf{k_G}$ is a matrix indexed in space, 
\begin{equation}
\mathbf{k_G}(\mathbf{x}, \mathbf{u};\mathbf{x}', \mathbf{u}')  : (\mathcal{S_X} \times \mathcal{S_U}) \times (\mathcal{S_X} \times \mathcal{S_U}) \rightarrow l \times l
\end{equation}
The GPs in the joint space are indexed by the design variables $\mathbf{x}$ and the sampled values of the random variables $\mathbf{U}$, written $\mathbf{u}$.  

The notations $F^{(t)}$ and $\mathbf{G}^{(t)}$ denote Gaussian processes conditioned by a set of $t$ observations: 
$\mathbf{f}^{(t)} = (f(\mathbf{x}^1,\mathbf{u}^1),\ldots, f(\mathbf{x}^t,\mathbf{u}^t))$ and  $\mathbf{g}^{(t)} = (\mathbf{g}(\mathbf{x}^1,\mathbf{u}^1),\ldots, \mathbf{g}(\mathbf{x}^t,\mathbf{u}^t))$ evaluated at $\mathcal D^{(t)} = \{(\mathbf{x}^i,\mathbf{u}^i)~,~i=1,..,t\}$. 

The conditional GP mean and covariance have the following expressions:
\begin{eqnarray}
m_F^{(t)}(\mathbf{x},\mathbf{u})) & = & m_F + \mathbf{k}_F(\mathbf{x},\mathbf{u})\mathbf{K_F}^{-1}(\mathbf{f}^{(t)}-\mathbf{m}_F) \\
k_F^{(t)}(\mathbf{x},\mathbf{u};\mathbf{x}',\mathbf{u}')   & = & k_F(\mathbf{x},\mathbf{u};\mathbf{x}',\mathbf{u}') - \mathbf{k}_F(\mathbf{x},\mathbf{u}) \mathbf{K_F}^{-1} \mathbf{k}_F^\top(\mathbf{x}',\mathbf{u}')
\end{eqnarray}
where $\mathbf{k}_F(\mathbf{x},\mathbf{u}) = \{ k_F(\mathbf{x},\mathbf{u},; \mathbf{x}^1,\mathbf{u}^1), \dots , k_F(\mathbf{x},\mathbf{u},; \mathbf{x}^t,\mathbf{u}^t)\}$,  $(\mathbf{K}_F)_{i,j} =  k_F(\mathbf{x}^i,\mathbf{u}^i; \mathbf{x}^j,\mathbf{u}^j)$ and $\mathbf{m}_F$ is a vector containing $t$ times the value $m_F$.

The same equations apply to the posterior mean and covariance of the GP $\mathbf{G}$: 
\begin{equation}
\mathbf{m}_{\mathbf{G}}^{(t)}(\mathbf{x},\mathbf{u}) = \mathbf{m}_\mathbf{G}(\mathbf{x}, \mathbf{u}) + \mathbf{K_G}(\mathbf{x},\mathbf{u}; \mathcal{D}^{(t)}) (\mathbf{K_G}( \mathcal{D}^{(t)}, \mathcal{D}^{(t)}) )^{-1} (\mathbf g^{(t)}-\mathbf{M}_\mathbf{G})
\end{equation}
\begin{multline}
\label{eq:K_G}
\mathbf{K_G}^{(t)}(\mathbf{x},\mathbf{u};\mathbf{x}',\mathbf{u}') = \mathbf{K_G}(\mathbf{x},\mathbf{u};\mathbf{x}',\mathbf{u}') \\
-  \mathbf{K_G}(\mathbf{x},\mathbf{u}; \mathcal{D}^{(t)}) (\mathbf{K_G}( \mathcal{D}^{(t)}, \mathcal{D}^{(t)}))^{-1} \mathbf{K_G}( \mathbf{x}',\mathbf{u}' ; \mathcal{D}^{(t)})^\top
\end{multline}
where $\mathbf{K_G}(\mathbf{x},\mathbf{u};\mathcal{D}^{(t)})$ is the $l \times tl$ matrix containing the covariance values over the $l$ outputs between the sample $\{\mathbf{x},\mathbf{u}\}$ and the training data set $\mathcal{D}^{(t)}$. 
Note that the dimensions of the vectors and matrices in the equations above are augmented with respect to the standard GP framework. 
For a prediction at a single point, $\mathbf{K_G}^{(t)}(\mathbf{x},\mathbf{u};\mathbf{x}',\mathbf{u}') \in \mathbb{R}^{l\times l}$, $\mathbf{K_G}(\mathcal{D}^{(t)},\mathcal{D}^{(t)}) \in \mathbb{R}^{tl\times tl}$, 
$\mathbf{K_G}( \mathbf{x},\mathbf{u};\mathcal{D}^{(t)}) \in   \mathbb{R}^{l\times tl}$.\\

Several approaches exist to model the structure of correlation between the $l$ components of $\mathbf{G}$. 
Without loss of generality
, the ``output as input'' (oai) approach is considered here. 
The multi-output GP is seen as a scalar GP defined on an increased input space, with an additional nominal input variable whose levels characterize one of the $l$ outputs. 
The additional nominal input is considered as an integer variable ranging from $1$ to $l$ and referred to with the letter $p$. 
Let $G_p$ and $G_{p'}$ be two components of $\mathbf{G}$, $Cov(G_p(\mathbf{x},\mathbf{u}),G_{p'}(\mathbf{x}',\mathbf{u}'))$ is assumed to be equal to $k^{oai}_{\mathbf{G}}(\{\mathbf{x},\mathbf{u},p\}, \{ \mathbf{x}',\mathbf{u}', p'\}) $ where
\begin{equation*}
 k^{oai}_{\mathbf{G}}(\{\mathbf{x},\mathbf{u},p\}, \{ \mathbf{x}',\mathbf{u}', p'\}) : (\mathcal{S_X} \times \mathcal{S_U} \times \mathbb{N}) \times (\mathcal{S_X} \times \mathcal{S_U} \times \mathbb{N}) \rightarrow \mathbb{R}
\end{equation*}
Note that the ``output as input'' covariance kernel $k^{oai}_{\mathbf{G}}(\cdot,\cdot)$ is a scalar function. Different ways of characterizing the aforementioned kernel exist. For instance,  $k^{oai}_{\mathbf{G}}(\cdot,\cdot)$ can be expressed as the product between the input dependent term and the output dependent term:
\begin{equation}
k^{oai}_{\mathbf{G}}(\{\mathbf{x}, \mathbf{u},p\}, \{ \mathbf{x}', \mathbf{u}', p'\}) = k_x(\mathbf{x}, \mathbf{u}; \mathbf{x}', \mathbf{u}') * k_p(p,p')
\label{eq:kernelTensorP}
\end{equation}

Although less common than the standard continuous kernels, a few kernels allowing to characterize the covariance function between the values of a discrete unordered variable such as $p$ exist in the literature. Different examples can be found in \cite{swiler2014surrogate, zhang2015computer, pelamatti2020overview}. Given that the main focus of this work is not to optimize the modeling performance of a multi-output GP model, but rather to study its effect within the framework of robust optimization, a single discrete kernel parameterization known as the hypersphere decomposition \cite{zhou2011simple} is considered for the remainder of the article. The underlying idea is detailed in Appendix \ref{sec:hypersphereKernel}.\\

Returning to the optimization Problem \eqref{newproblem}, models for the functions $z$ and $c$ in the design space $\mathcal{S_X}$ are deduced from the models of the objective function and the constraints in the joint space. 
Since the expectation is a linear operator, the averaged process over the uncertain domain remains Gaussian. 
A GP modeling the average objective function has therefore been introduced:
\begin{equation*}
Z^{(t)}(\mathbf{x}) =\mathbb E_{\mathbf{U}}[F^{(t)}{(\mathbf{x},\mathbf{U})}]
\end{equation*}
The mean and covariance function of $Z$ are:
\begin{eqnarray}
m^{(t)}_Z(\mathbf{x}) & = & \int_{\mathbb{R}^m} m^{(t)}_F(\mathbf{x},\mathbf{u}) \rho_{\mathbf{U}}(\mathbf{u}) d\mathbf{u}, \\
k^{(t)}_Z(\mathbf{x},\mathbf{x'}) &=& \iint \limits_{\mathbb{R}^m} k^{(t)}_F(\mathbf{x},\mathbf{u},\mathbf{x'},\mathbf{u'}) \rho_{\mathbf{U}}(\mathbf{u'}) \rho_{\mathbf{U}}(\mathbf{u'}) d\mathbf{u} d\mathbf{u'}
\label{ZGP}
\end{eqnarray}
The integrals in \eqref{ZGP} can be evaluated analytically under some assumptions (see \cite{janusevskis_jogo_2012}). 
When this is not feasible a quadrature scheme can be used to approximate these integrals. 

The probability of feasibility with a given safety margin $\alpha$ is estimated with the GPs $\mathbf{G}$ through the process
\begin{equation*}
C^{(t)}(\mathbf{x}) = 1 - \alpha - \mathbb E_\mathbf{U}[\indic_{\mathbf{G} (\mathbf{x}, \mathbf{U}) \leq \mathbf{0} } ]
\end{equation*}
where $ \mathbf{0} $ is a $l \times 1$ vector of zeros.  
It is important to note that, by construction, $C$ is not Gaussian. 
This characteristic makes the analytical tractability of the proposed method more complicated, as it will be discussed later.

\subsection{Robust Bayesian optimization}
\label{RobBayOpt}
Let us now turn to the adaptation of Bayesian optimization to Problem \eqref{newproblem}.
We generalize the approach described in \cite{elamri2021sampling} to multi-output 
models of the constraints.
The algorithm performs both optimization and sampling in that the design variables and the uncertain variables at which the exact problem functions are to be evaluated are determined in the same iteration. The working principle of this algorithm is to maximize the feasible improvement, which stands as the acquisition criterion, but in a one-step-ahead uncertainty reduction strategy. 
The method is detailed in the following paragraphs.
\\
The most promising location of the design space, denoted $x_{targ}$ for targeted $x$, is the set of design variables which maximizes the Expected Feasible Improvement (EFI) : 
\begin{equation}\label{eq:EFI2}
\mathbf{x}_{targ} = \mbox{arg} \min_{\mathbf{x}} \mathbb E[FI^{(t)}(\mathbf{x})]
\end{equation}
where $FI^{(t)}(\mathbf{x}) = \left(z_{min}^{feas} - Z^{(t)}(\mathbf{x})\right)^+ \indic_{\{  C^{(t)}(\mathbf{x}) \leq 0   \}}$. \\
By exploiting the independence between the objective function and the constraints, Equation \eqref{eq:EFI2} can be written as:
\begin{eqnarray}\label{eq:EFI3}
\mathbf{x}_{targ} = \mbox{arg} \min_{\mathbf{x}} EI^{(t)}(\mathbf{x}) \mathbb P(C^{(t)} (\mathbf{x}) \leq 0 ). 
\end{eqnarray}\\
The first term of \eqref{eq:EFI3}, the Expected Improvement, is a common acquisition criterion and it is written
\begin{equation}
EI^{(t)}(\mathbf{x}) = (z_{\min}^\text{feas} - m_Z^{(t)}(\mathbf{x})) \Phi\bigg(\frac{z_{\min}^\text{feas} - m_Z^{(t)}(\mathbf{x})}{\sigma_Z^{(t)}(\mathbf{x})}\bigg) + \sigma_Z^{(t)}(\mathbf{x}) \phi\bigg(\frac{z_{\min}^\text{feas} - m_Z^{(t)}(\mathbf{x})}{\sigma_Z^{(t)}(\mathbf{x})}\bigg),
\label{eq:EI}
\end{equation}
where $\sigma_Z^{(t)}(\mathbf{x})$ is the standard deviation of the conditioned GP $Z$ at $\mathbf{x}$, while $\Phi$ and $\phi$ are the normal cumulative distribution and density functions, respectively.  
$z_{\min}^\text{feas}$ designates the best feasible target value for the mean function. 

An expression of Equation \eqref{eq:EFI3} is given in \cite{elamri2021sampling} but it must be adapted to the multi-output context. 
The second term of \eqref{eq:EFI3}, $\mathbb{P}(C^{(t)}(\mathbf{x}) \leq 0)$,  was estimated numerically with the help of a Monte Carlo routine. 
This computation is extended to the multi-output model: the simulation of independant GPs is replaced by the simulation of multi-output GP trajectories. 
In Equation \eqref{eq:EI}, the expression of $z_{\min}^\text{feas}$ involves the cumulative density function (cdf) of a multivariate Gaussian distriubution instead of 
the product of $l$ Gaussian cdf's.
The details of these calculations are provided in Appendix \ref{sec:detailsFeasImprov}. \\ 
 
Once $x_{targ}$ is determined, the design of experiments $\mathcal{D}^{(t)}$ is enriched with the next iterate $\{(\mathbf{x}^{t+1},\mathbf{u}^{t+1})\}$.
This iterate is such that a proxy $S$ of the one-step-ahead variance of the acquisition criterion (the feasible improvement) is minimized at $x_{targ}$. $S$ will serve as a Sampling criterion for the random variables, hence the associated symbol.
In \cite{janusevskis_jogo_2012}, it is shown that the optimal values of $\mathbf{x}^{t+1}$ tend to be very close to $x_{targ}$. 
For this reason, in order to reduce the computational overhead, it is assumed that $\mathbf{x}^{t+1} = x_{targ}$. 
Then, only $\mathbf{u}^{t+1}$ needs to be computed, which is done as follows:
\begin{equation}
\label{eq:uplus1}
\mathbf{u}^{t+1} = \arg \min_\mathbf{\tilde{u}} S(\mathbf{x}_{targ},\mathbf{\tilde{u}}) ~.
\end{equation}

A global view of the entanglement of the above steps within a constrained, robust, optimization method can be gained by looking at Algorithm \ref{alg:template} for which other options will be introduced later (see Section \ref{sec:theAlgos}).\\

The proxy to the one-step-ahead variance of the feasible improvement is written
\begin{multline} 
\label{eq:proxy}
S(\tilde{\mathbf{x}},\mathbf{\tilde{u}}) =~Var(\big(z_{\min}^{\mbox{feas}} - Z^{(t+1)}(\mathbf{x}_{targ})\big)^+) \\
\int_{\mathbb{R}^m} Var\big( \indic_{\{\mathbf{G}^{(t+1)}(\mathbf{x}_{targ},\mathbf{u}) \leq 0 \}} \big) \rho_{\mathbf{U}}(\mathbf{u}) d\mathbf{u} 
\end{multline}
where $Z^{(t+1)}$ and $\mathbf{G}^{(t+1)}$ are conditioned on $\mathcal{D}^{(t)} \cup \{(\tilde{\mathbf{x}},\mathbf{\tilde{u}})\}$.\\

The first term of $S(\tilde{\mathbf{x}},\mathbf{\tilde{u}})$, i.e. $Var\left( \big(z_{\min}^{\mbox{feas}} - Z^{(t+1)}(\mathbf x_{targ})\big)^+\right)$, corresponds to the improvement variance. 
Its expression has some resemblance to the one of the expected improvement and is given in \cite{elamri2021sampling}. 
It is expressed with the probability and cumulated density functions of the standard Gaussian distribution (see Appendix \ref{sec:VarItp1calc}).

In order to compute the second term of $S(\tilde{\mathbf{x}},\mathbf{\tilde{u}})$ with a multi-output constraint model, one needs to determine $\mathbf{G}^{(t+1)}$ by conditioning $\mathbf{G}^{(t)}$ on the new sample $\{\tilde{\mathbf{x}},\tilde{\mathbf{u}}\}$ while taking into account the correlation between the components of $\mathbf{G}^{(t)}$. 
By drawing on the Kriging Believer assumption\footnote{
The Kriging Believer assumption states that the value of a function at an unobserved point is equal to the kriging prediction at that point. For our constraints, it means $g_p(\tilde{\mathbf x},\tilde{\mathbf u}) = {{\mathbf m}_{\mathbf G}^{(t)}}_p(\tilde{\mathbf x},\tilde{\mathbf u}) ~,~p=1,\ldots,l$.}, the updated prediction is 
\begin{eqnarray}
\mathbf{m_G}^{(t+1)}(\mathbf{x}_{targ},\mathbf{u})  &=& \mathbf{m_G}^{(t)}(\mathbf{x}_{targ},\mathbf{u}) \nonumber \\ 
\mathbf{K_G}^{(t+1)}(\mathbf{x}_{targ},\mathbf{u} ; \mathbf{x}_{targ},\mathbf{u}) &=  & \mathbf{K_G}^{(t)}(\mathbf{x}_{targ},\mathbf{u};\mathbf{x}_{targ},\mathbf{u}) \label{MOUpdate} \\
&&-\mathbf{K_G}^{(t)}(\mathbf{x}_{targ},\mathbf{u};\tilde{\mathbf{x}},\mathbf{\tilde{u}} ) 
\mathbf{K}_{new}^{-1}
 \mathbf{K_G}^{(t)}(\mathbf{x}_{targ},\mathbf{u};\tilde{\mathbf{x}},\mathbf{\tilde{u}} )^\top \nonumber
\end{eqnarray}
with $\mathbf{K}_{new} =  \mathbf{K_G}^{(t)}(\tilde{\mathbf{x}},\mathbf{\tilde{u}};\tilde{\mathbf{x}},\mathbf{\tilde{u}})$ the $l \times l$ matrix containing the covariance of the vector $\mathbf G(\tilde{\mathbf{x}},\mathbf{\tilde{u}})$, 
i.e., the covariance values between the different outputs of the GP evaluated at $(\tilde{\mathbf{x}},\mathbf{\tilde{u}})$.

Now, $\indic_{\mathbf{G}^{(t+1)}(\mathbf{x}_{targ}, \mathbf{u}) \leq \mathbf{0} }$ is a Bernoulli variable with parameter $p(\mathbf{u})$ equal to 
\begin{equation}
p(\mathbf{u}) = \Phi\left( \mathbf{0}-\mathbf{m_G}^{(t)}(\mathbf{x}_{targ},\mathbf{u}),\mathbf{K}^{(t+1)}_{\mathbf{G}}(\mathbf{x}_{targ},\mathbf{u};\mathbf{x}_{targ},\mathbf{u}) \right) ~,
\end{equation}
where $\Phi(\cdot)$ is the cumulative distribution function of a multivariate Gaussian distribution.
As a consequence, 
\begin{eqnarray}
\int_{\mathbb{R}^m} Var (\indic_{\mathbf{G} ^{(t+1)}(\mathbf{x}_{targ}, \mathbf{u}) \leq \mathbf{0} }) \rho (\mathbf{u}) d\mathbf{u} = \int_{\mathbb{R}^m} p(\mathbf{u})(1-p(\mathbf{u}))  \rho_\mathbf{U} (\mathbf{u}) d\mathbf{u} 
\end{eqnarray}\\

As it will be seen in Section \ref{sec:illustration}, modeling the correlation between the constraints increases the accuracy of the surrogate. 
However, the difference in performance between the two variants (with or without constraint correlation) will largely depend on the actual correlation value between the modeled functions. 
If the constraints are independent for example, the multi-output model will be more complex than the standard approach while providing no additional benefit. 
More generally, the potential gain in performance of the coupled model depends on the ability of the correlation structure to describe the constraints relationship.
In any case, accounting for constraints couplings results in a substantial increase of the covariance matrix size, which may become a practical bottleneck with large data sets. When this happens, a sparse GP covariance structure might be indicated.

\vspace{12pt}

\subsection{Illustration of feasibility estimation with a multi-output GP}
\label{sec:illustration}

In order to illustrate the behaviour of the multi-output model of the constraints, an analytical example where the constraints are designed to be linked is considered:
\begin{equation}
\begin{split}
& \min_{\mathbf{x}} \mathbb E_U[f(\mathbf{x},\mathbf{U})] ~\text{ such that }~ \mathbb P(g_p(\mathbf{x},\mathbf{U}) \leq 0, p = 1,2) \geq 0.95 \\
& \text{where } \\
& f(x_1,x_2,u_1,u_2)  =  5(x_1^2+x_2^2) - (u_1^2+u_2^2) + x_1(u_2-u_1+5) + x_2(u_1-u_2+3)  \\ 
& g_1(x_1,x_2,u_1,u_2)  =   -x_1^2 + 5x_2 - u_1 + u_2^2 -1 \\ 
& g_2(x_1,x_2,u_1,u_2)  =  g_1(x_1,x_2,u_1,u_2)(x_1+5)/5 - u_1 - 1
\end{split}
\label{eq:4d_analytic}
\end{equation}
with $\mathbf{x} \in [-5,5]^2$ and $\mathbf{U} \sim \mathcal{U}([-5,5]^2)$. 
A correlation is induced between the two constraints by making the expression of the second constraint a transformation of the expression of the first constraint.

\begin{figure}[h!]
\begin{minipage}{0,5\textwidth}
\includegraphics[width=\textwidth]{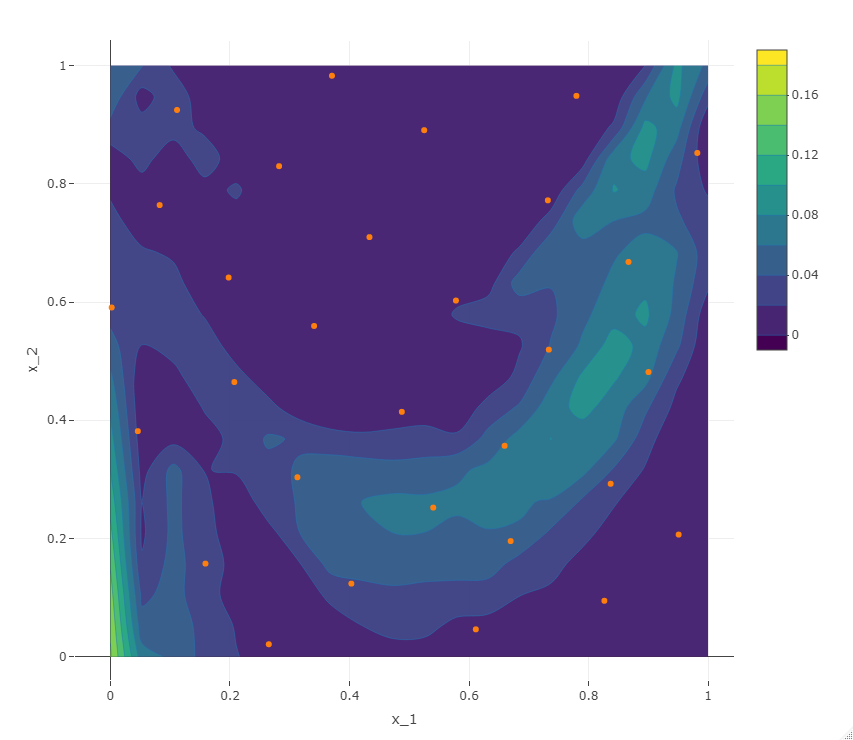}
\end{minipage} %
\begin{minipage}{0,5\textwidth}
\includegraphics[width=\textwidth]{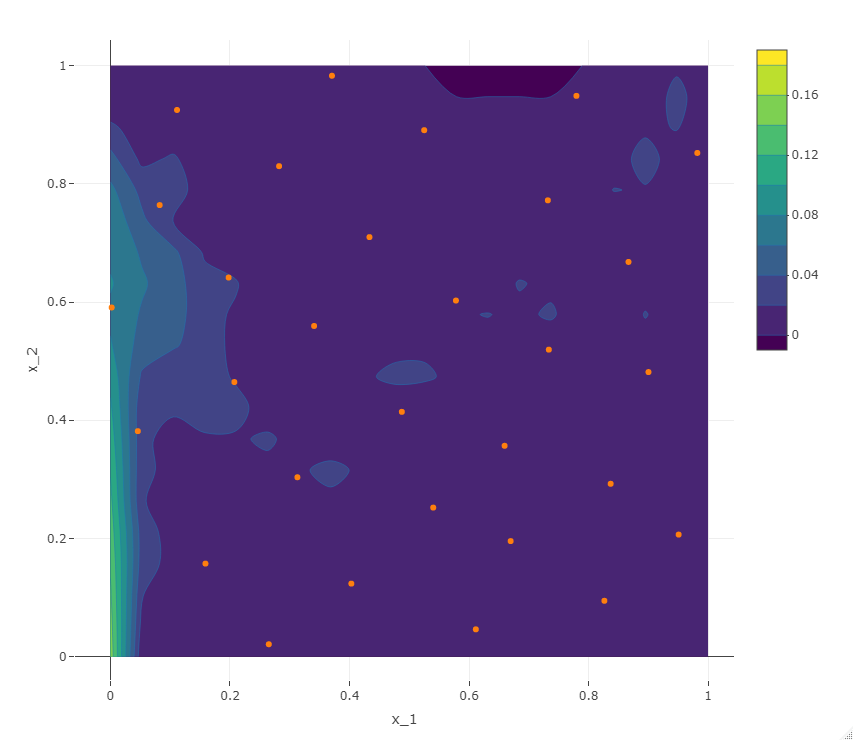}
\end{minipage}
\caption{Contour plots of the error in the probability of feasibility on the Problem \eqref{eq:4d_analytic}. 
Left: independent modeling of each constraint. Right: multi-output GP model of dependent constraints. 
The training points in the $\mathcal{S_X}$ space are indicated with the orange dots.
}
\label{PoFMOMod}
\end{figure}

\begin{figure}[h!]
\centering
\includegraphics[width=0.7\textwidth]{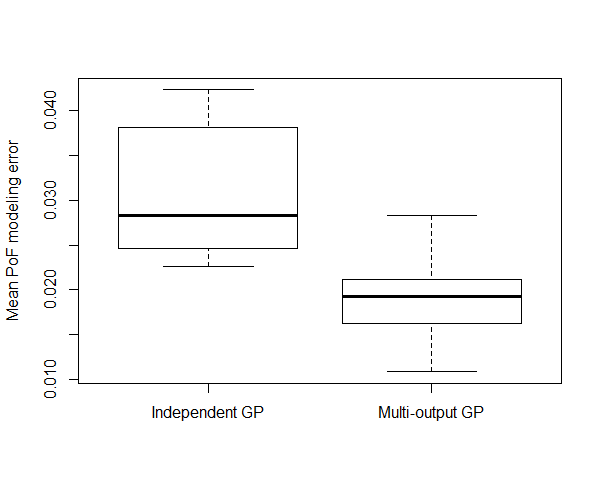}
\caption{Average error in probability of feasibility on Problem~\eqref{eq:4d_analytic}. 
The left boxplot is obtained with a an independent GP for each constraint, whereas the right boxplot corresponds to a correlated multi-output GP.
}
\label{BoxPlotPoFError}
\end{figure}

Depending on whether an independent or a dependent model of the constraints is considered, the following expression of the Probability of Feasibility (PoF) stands:
\begin{eqnarray*}
\mbox{ Independent constraints: } \ \ \text{PoF}(\mathbf{x}) & = & E_\mathbf{U}[\indic_{\cap_{p = 1}^2 \{G_p (\mathbf{x}, \mathbf{U}) \leq 0 \} }] \\
\mbox{ Dependent constraints: } \ \  \text{PoF}(\mathbf{x}) & = & E_\mathbf{U}[\indic_{\mathbf{G} (\mathbf{x}, \mathbf{U}) \leq \mathbf{0} } ]
\end{eqnarray*}
The modeling performance is calculated as the absolute difference between the PoF computed with the GP models, and the one computed on the true function at $N=1000$ $\mathbf{u}^j$'s sampled with respect to $\mathcal{U}([-5,5]^2)$,
\begin{equation*}
\text{PoF error}(\mathbf x) ~=~ |\text{PoF}(\mathbf{x}) - \frac{1}{N}  \sum_{j = 1}^{N} 1_{\cap_{p = 1}^2 \{g_p (\mathbf{x}, \mathbf{u}^j) \leq 0 \} } | ~.
\end{equation*} 
The results are plotted in Figure~\ref{PoFMOMod} where, for the sake of clarity, the input variables are normalized between 0 and 1 in the plots.
A data set of 30 samples per constraint function is used to train the scalar and the multi-output GPs.
The Figure shows that the areas of $\mathcal{S_X}$ characterized by the largest accuracy (lowest error) are the same in both models. 
Overall however, the multi-output GP provides a lower PoF estimation error than the independent model. 
The gain in accuracy is particularly clear in the parabolic region where the constraints are near 0, which is the considered feasibility threshold.
The advantage of the multi-output constraint model is further highlighted in Figure \ref{BoxPlotPoFError} where the mean PoF error computed on 400 randomly sampled values in $S_X$ over 10 repetitions of the DoE and GPs training is shown for both approaches.
The difference in modeling performance can be explained by the fact that the multi-output GP benefits from a larger data set since it processes information from the training data of both constraints. 
In the context of optimization, it is reasonable to assume that a more accurate feasibility prediction will allow a better search for robust optimal solutions, which highlights the potential benefit of a multi-output model of the constraints.

\section{Extensions to robust Bayesian optimization enabled by constraints coupling}
\label{ProxExt}
Two extensions to robust Bayesian optimization are now proposed. 
The first extension is the possibility of refining the objective and constraint functions with different sets of random variables (but identical design variables point). 
The second extension allows to select a single constraint to evaluate at $x_{targ}$ in order to save constraints evaluations. 
The chosen constraint provides the largest amount of information regarding the feasibility at the candidate optimum location $x_{targ}$. 
These extensions are relevant when each output, at least the objectif function and the set of constraints, can be obtained through separate simulators
in which case a goal oriented sampling strategy that sparingly selects which simulator should be called will save computations.

\subsection{Independent selection of uncertain variables}
The robust Bayesian optimization algorithm of Section \ref{Multioutput} associates to each target design variable, $\mathbf{x}_{target}$, a single vector of random variables $\mathbf{u}^{t+1}$ that will be used to refine the GPs of both the objective and the constraints functions.
In most cases, however, the most informative $\mathbf u$'s from the objective function improvement and from the constraint satisfaction perspectives are different. 
Selecting a single $\mathbf{u} \in \mathbb{R}^m$ for a better accuracy in, both, the improvement and the failure probability, is necessarily a compromise. 
For this reason, we introduce the possibility of selecting different vectors, $\mathbf{u}_f$ and $\mathbf{u}_g$, for the objective and the constraint functions. 
Again, the underlying assumption in this duplication of the random samples is that the computationally intensive codes used to evaluate the objective and the constraints can be called independently. \\

Remember that the acquisition function of Equation \eqref{eq:proxy},
\begin{equation*}
S(\mathbf{x}_{targ},\mathbf{u}) = Var((z_{min}^{feas}-Z^{(t+1)}(\mathbf{x}_{targ}))\int_{\mathbb{R}^m} Var (1_{\mathbf{G}^{(t+1)} (\mathbf{x}_{targ}, \mathbf{u}) \leq \mathbf{0} }) \rho (\mathbf{u}) d\mathbf{u} ~,
\end{equation*}
is composed of two terms, the first related to the objective function and the second to the constraints. 
It can be decomposed into two independent terms, $S_f$ and $S_g$ 
involving different samples of the random variables : 
\begin{eqnarray}
S_f(\mathbf{x}_{targ},\mathbf{u}_f) & =  & Var((z_{min}^{feas}-Z^{(t+1)}(\mathbf{x}_{targ})^+)   \label{eq:S1} \\\nonumber
 & =  & Var((z_{min}^{feas}-Z(\mathbf{x}_{targ}))^+| F(\mathcal{D}^{(t)}_f) = f^{(t)}, F(\mathbf{x}_{targ},\mathbf{u}_f)  )\\
S_g(\mathbf{x}_{targ},\mathbf{u}_g) & = & \int_{\mathbb{R}^m} Var (1_{\mathbf{G} ^{(t+1)}(\mathbf{x}_{targ}, \mathbf{u}) \leq \mathbf{0} }) \rho (\mathbf{u}) d\mathbf{u} \label{eq:S2} \\
& = & \int_{\mathbb{R}^m} Var (1_{\mathbf{G} (\mathbf{x}_{targ}, \mathbf{u})\leq \mathbf{0} | \mathbf{G}(\mathcal{D}^{(t)}_g) = \mathbf{g}^{(t)}, \mathbf{G}(\mathbf{x}_{targ},\mathbf{u}_g) }) \rho (\mathbf{u}) d\mathbf{u} \nonumber
\end{eqnarray}
where $\mathcal{D}^{(t)}_f$ and $\mathcal{D}^{(t)}_g$ are design of experiments after $t$ evaluations that contain different samples of $u$.
More precisely if $t_{init}$ denotes the size of the initial DOE, we have $\mathcal{D}^{(t_{init})}_f=\mathcal{D}^{(t_{init})}_g =\mathcal{D}^{(t_{init})}$. For $t>t_{init}$ , the DOE for $f$ (for $g$, respectively) evolves as  $\mathcal{D}^{(t)}_f=\mathcal{D}^{(t-1)}_f \cup \{(\mathbf x^t,\mathbf u^t_f)\}$,  ($\mathcal{D}^{(t)}_g=\mathcal{D}^{(t-1)}_g \cup \{(\mathbf x^t,\mathbf u^t_g)\}$, respectively). \\

As a consequence, it is possible to separately optimize $S_f$ and $S_g$, thus obtaining two separate vectors of uncertain variables, one for the next objective function and one for the next constraints evaluation:
\begin{eqnarray}
\mathbf{u}_f^{t+1} & = & \arg\min_{\mathbf{u}_f} S_f(\mathbf{x}_{targ},\mathbf{u}_f) \label{eq:ufnext} \\
\mathbf{u}_g^{t+1} & = & \arg\min_{\mathbf{u}_g} S_g(\mathbf{x}_{targ},\mathbf{u}_g) \label{eq:ugnext}
\end{eqnarray}

This enables the algorithm to more efficiently refine the surrogate models of the problem functions : the set of variables $\mathbf{u}_f$   minimizes the one-step-ahead variance of the improvement, and  $\mathbf{u}_g$  minimizes the one-step-ahead variance of the probability of feasibility.
The generic Algorithm \ref{alg:template} may help the reader in positioning the constraint selection operation within the global optimization procedure (see Section \ref{sec:theAlgos}).

\begin{figure}[h!]
\centering
\includegraphics[width=0.6\textwidth]{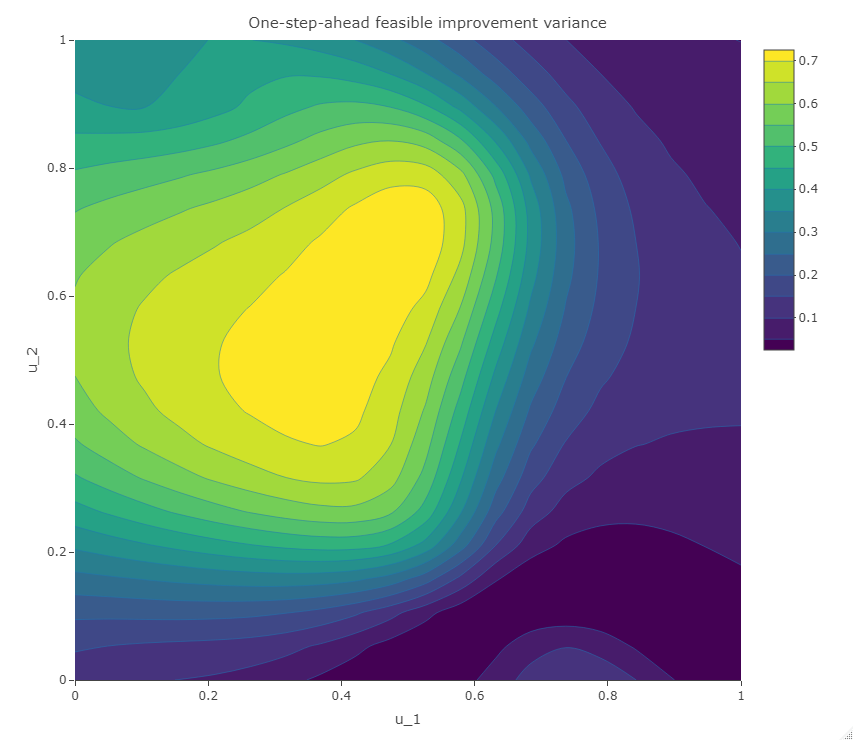}
\caption{Contour plot of $S$, the proxy of the one-step-ahead feasible improvement variance computed at a given $\mathbf{x}_{target}$.  
\label{FeasImp} }
\end{figure}

\begin{figure}[h!]
\begin{minipage}{0,5\textwidth}
\includegraphics[width=\textwidth]{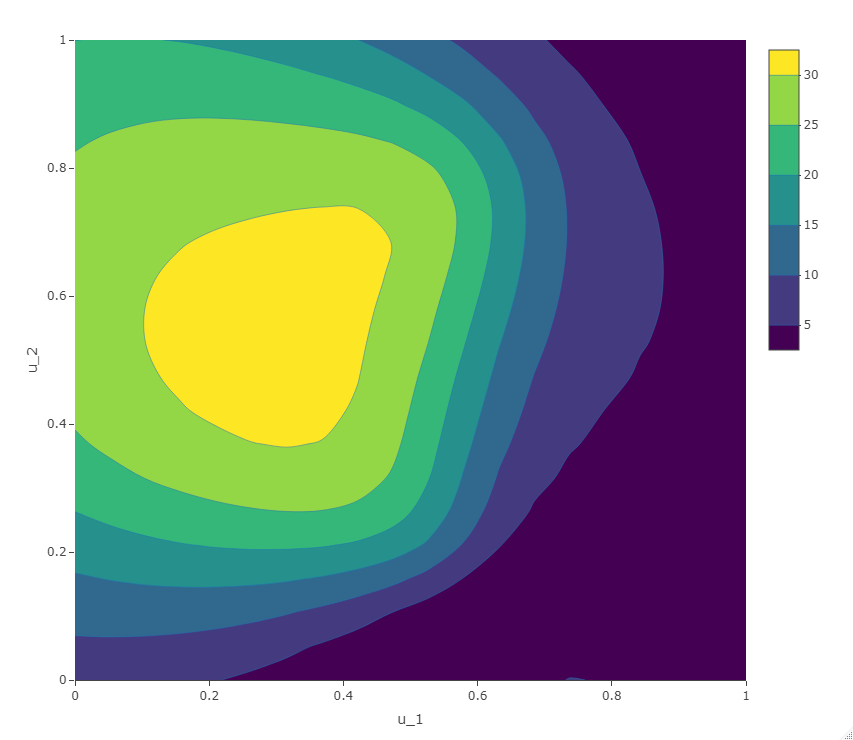}
\end{minipage} %
\begin{minipage}{0,5\textwidth}
\includegraphics[width=\textwidth]{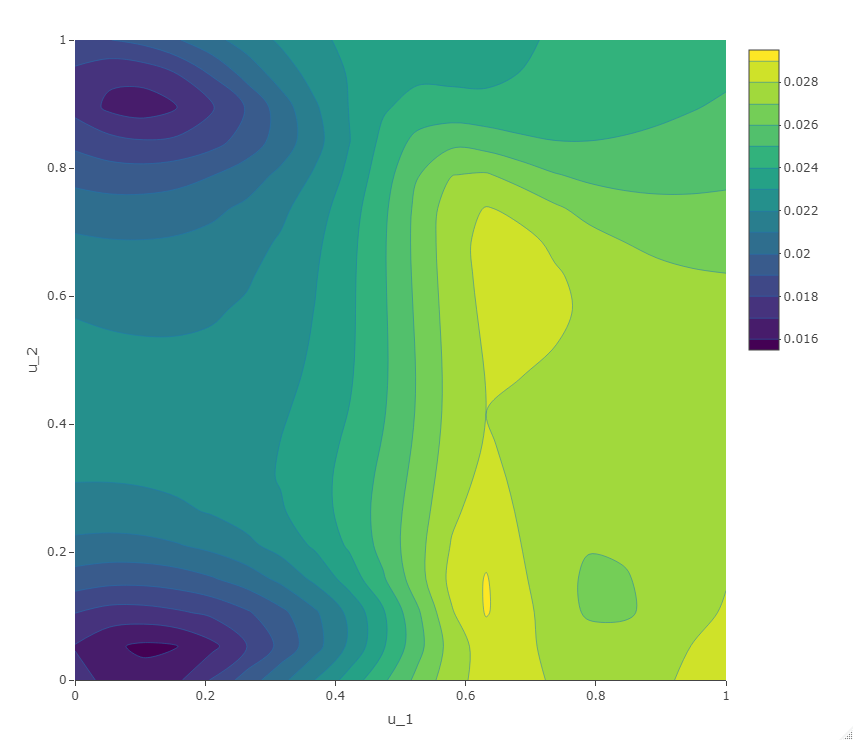}
\end{minipage} 
\caption{Contour plots of the one-step-ahead variance of the improvement (left) and the feasibility probability (right), computed at a given $\mathbf{x}_{target}$.
\label{SC1SC2}
}
\end{figure}

\vskip\baselineskip
To highlight the effect of this modification, we return to the analytical Problem \eqref{eq:4d_analytic}.
A data set of size $t=20$ is created, GPs of the objective function and the two constraints are learned, and the most promising location in the design space $\mathbf{x}_{target}$ determined by maximization of the Expected Feasible Improvement (Equation \eqref{eq:EFI3}).
The GPs and $\mathbf{x}_{target}$ can then be used in the formulas for either $S$, the proxy to the one-step ahead feasible improvement variance, or its terms $S_f$ and $S_g$. At that point, these sampling criteria only depend on the random parameters $u_1$ and $u_2$. 
A contour plot of $S$ as a function of $u_1$ and $u_2$ is drawn in Figure \ref{FeasImp}.
It is meant to be compared to the contours of $S_f$ and $S_g$, the one-step-ahead variance of the expected improvement and of the feasibility probability, given in Figure \ref{SC1SC2}.

Samples of the random parameters are generated by minimization of the one-step-ahead variances, $S$, $S_f$ or $S_g$.
When comparing the feasible improvement variance proxy, $S$, with its two factors, $S_f$ and $S_g$, 
it is observed that the improvement variance ($S_f$) is the dominating factor in $S$ and drives the selection of the next samples 
$\mathbf{u}^{t+1}$ in the generic algorithm of Equation \eqref{eq:uplus1}.
However, such a $\mathbf{u}^{t+1}$ will poorly reduce the variance of the probability of feasibility ($S_g$ map).
In the above example, $\mathbf{u}^{t+1} = \{1,0\}$ when the best sample for the constraints is $\mathbf{u}^{t+1}_g = \{0.1,0.05\}$.

Therefore a unique $\mathbf u$ may not be the best strategy to efficiently learn where is the feasible domain.

\subsection{Iterative optimal constraint selection}
\label{ConstSel}

When the constraints relate to different physics or disciplines, 
it is common to have separate simulation codes for each constraints function (or group of constraints).

A direct extension of the previous section would then be to consider different instances of the random variables for each of the constraint. The acquisition function could be detailed for each of the constraint random sample $\mathbf{u}_{g_i}$

\begin{eqnarray}
&& S_g^{all}(\mathbf{x}_{targ},\mathbf{u}_{g_1}, ..., \mathbf{u}_{g_l}) =   \\ &&
\int_{\mathbb{R}^m} Var (1_{\mathbf{G} (\mathbf{x}_{targ}, \mathbf{u})\leq \mathbf{0}} | \mathbf{G}(\mathcal{D}^{(t)}_g) = g^{(t)}, G_{1}(\mathbf{x}_{targ},\mathbf{u}_{g_1}), \ldots, G_{l}(\mathbf{x}_{targ},\mathbf{u}_{g_l})) \rho (\mathbf{u}) d\mathbf{u} \nonumber 
\end{eqnarray}
where $\mathcal{D}^{(t_{init})}_g=\{(\mathbf x^i,\mathbf u^i,1),\ldots,(\mathbf x^i,\mathbf u^i,l), \forall 1\leq i \leq t_{init}\}$ 
and for $t>t_{init}$, $\mathcal{D}^{(t)}_g=\mathcal{D}^{(t-1)}_g \cup \{(\mathbf x^t,\mathbf u^t_{g_1},1),\ldots,(\mathbf x^t,\mathbf u^t_{g_l},l) \}$.
Remark that the 
covariance matrix to be inverted in the update Equation \eqref{MOUpdate}, $\mathbf{K}_{new}$, is slightly more complex as it is computed with respect to data samples characterized by different values of both $\mathbf{u}$ and $p$.\\

When the codes of each constraint are autonomous, a more appealing extension for numerical efficiency is to select only one constraint to update at each iteration. Indeed,
the refinement of some constraints may not provide any additional information when the surrogate models are already accurate enough to determine feasibility at a given location.  
Furthermore, the step-by-step selection of the most relevant constraint to evaluate prevents an unnecessary increase of the covariance matrix size. 
This is important in Bayesian optimization because the dominating numerical cost comes from the repeated $\mathcal O(t^3)$ inversion of the covariance matrix during the likelihood maximizations.
Adding constraint evaluations one at a time is particularly tempting when the constraints are modeled with a multi-output GP since all the GPs improve with that one evaluation.

In order to select the constraint, the one-step-ahead variance of the probability of feasibitily, $S_g$ (see Equation \eqref{eq:S2}), 
is modified by including the index $p$ of the constraint function considered, $g_p$, 
and the associated sampling location  $\mathbf{u}_{g}$  :

\begin{equation}
S_g(\mathbf{x}_{targ},\mathbf{u}_{g},p)   =  \int_{\mathbb{R}^m} Var (1_{\mathbf{G} (\mathbf{x}_{targ}, \mathbf{u})\leq \mathbf{0} | \mathbf{G}(\mathcal{D}^{(t)}_g) = g^{(t)}, G_p(\mathbf{x}_{targ},\mathbf{u}_{g}) }) \rho_\mathbf{U} (\mathbf{u}) d\mathbf{u}
\label{SC2_ConstSel}
\end{equation}
where for $t>t_{init}$, $\mathcal{D}^{(t)}_g=\mathcal{D}^{(t-1)}_g \cup \{(\mathbf x^t,\mathbf u^t_{g},p^t)\}$.
In the above Equation, the one-step-ahead update of the GP includes a single additional data sample for the constraint designated by $p$. 
In this case, the covariance matrix of $\mathbf G(\mathcal D^{(t)})$ has a size smaller than $tl\times tl$ 
and its update coming from Equation \eqref{MOUpdate} is 
\begin{equation*}
\mathbf{K_G}^{(t+1)}(\mathbf{x}_{targ},\mathbf{u} ; \mathbf{x}_{targ},\mathbf{u})  =   \mathbf{K_G}^{(t)}(\mathbf{x}_{targ},\mathbf{u}; \mathbf{x}_{targ},\mathbf{u})  - \mathbf{w} k_{new}^{-1}  \mathbf{w}^\top 
\end{equation*}
with $k_{new}$ a scalar and $\mathbf{w}$ a $l \times 1$ vector,
\begin{eqnarray*}
k_{new} &=& Cov(G_p^{(t)}(\mathbf{x}_{targ},\mathbf{u}_{g}),G_p^{(t)}(\mathbf{x}_{targ},\mathbf{u}_{g})) \\ 
\mathbf{w}& = & \left( Cov(G_1^{(t)}(\mathbf{x}_{targ},\mathbf{u}),G_p^{(t)}(\mathbf{x}_{targ},\mathbf{u}_{g}) ), \ldots, Cov(G_l^{(t)}(\mathbf{x}_{targ},\mathbf{u}),G_p^{(t)}(\mathbf{x}_{targ},\mathbf{u}_{g}))\right)^\top \\ 
\end{eqnarray*}

For selecting the constraint at update time, the one-step-ahead feasibility variance is the acquisition function which is simultaneously optimized with respect to $\mathbf{u}_{g}$ and $p$.
The solution is a pair characterizing both the constraint and the parameter in the $\mathcal U$ space which provides the largest amount of information regarding the feasibility of $\mathbf x_{targ}$,
\begin{equation}
p^{t+1}, \mathbf{u}_{g}^{t+1} = \arg\min_{p,\mathbf{u}_{g}} S_g(\mathbf{x}_{targ},\mathbf{u}_{g},p) \quad.
\label{eq:minS2oat}
\end{equation}
At each iteration, the random sample of the objective function, $\mathbf{u}_f^{t+1}$ comes from the one-step-ahead improvement variance minimization as expressed in Equation \eqref{eq:ufnext}.
An evaluation of the objective function is performed together with each evaluation of the selected constraint at $(\mathbf{x}_{targ},\mathbf{u}_f^{t+1})$ and $(\mathbf{x}_{targ},\mathbf{u}_{g}^{t+1})$, respectively. 
Hence, for each objective function evaluation there is a single constraint function evaluation, rather than $l$ distinct constraint function evaluations. A reason for this is that once it is known that a constraint is satisfied or violated, there is no use in fine tuning it. On the contrary, further adjusting the objective function model in high-performance regions is usually useful. 
\vspace{12pt}

\begin{figure}[h!]
\begin{minipage}{0,5\textwidth}
\includegraphics[width=\textwidth]{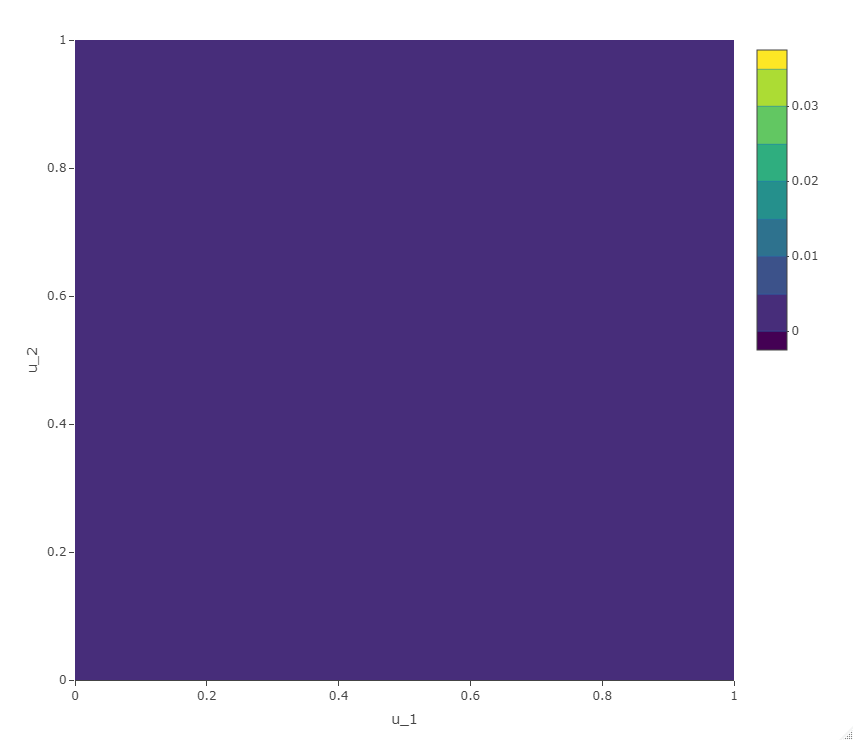}
\end{minipage} %
\begin{minipage}{0,5\textwidth}
\includegraphics[width=\textwidth]{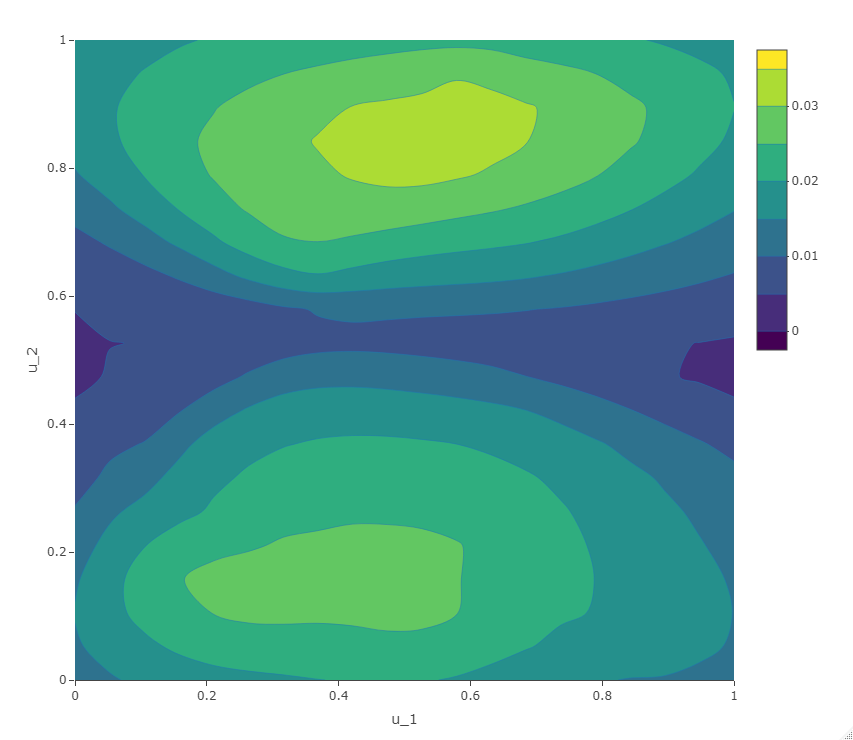}
\end{minipage} 
\caption{One-step-ahead feasibility variance reduction which would be obtained by refining the model of either $g_1$ (left) or $g_2$ (right) at the value $(\mathbf x_{targ},u_1,u_2)$, Problem \eqref{eq:4d_analytic}}
\label{VarRedg}
\end{figure}

In order to highlight the potential benefit of selecting the most useful constraint at every iteration, 
the analytical test case described in Equations \eqref{eq:4d_analytic} is considered anew. 
For illustrative purposes, the $\mathbf{x}_{target}$ value ($\{0,0.3\}$) is manually selected so that the constraint $g_1$ is almost certainly satisfied, 
whereas the probability of feasibility for $g_2$ is approximately $0.6$. Each constraint is modeled by an independent GP, trained by relying on a data set of 30 randomly generated samples.
The quantity of interest is the reduction in one-step-ahead variance of the feasibility which would be obtained by refining the model of either $g_1$ (when $p=1$) or $g_2$ ($p=2$),
\begin{multline*}
 \int_{\mathbb{R}^m} Var (1_{\mathbf{G} (\mathbf{x}_{targ}, \mathbf{u})\leq \mathbf{0} | \mathbf{G}(\mathcal{D}^{(t)}_g) = g^{(t)}}) \rho (\mathbf{u}) d\mathbf{u} ~-~ \\
\int_{\mathbb{R}^m} Var (1_{\mathbf{G} (\mathbf{x}_{targ}, \mathbf{u})\leq \mathbf{0} | \mathbf{G}(\mathcal{D}^{(t)}_g) = g^{(t)}, \mathbf{G}(\mathbf{x}_{targ},\mathbf{u}_g, p) }) \rho (\mathbf{u}) d\mathbf{u} ~.
\end{multline*}
This quantity is plotted for both values of $p$ (i.e., for either a $g_1$ or a $g_2$ enrichment) in Figure \ref{VarRedg}.
The left part of the Figure shows that refining the $g_1$ constraint which is surely feasible at $\mathbf{x}_{targ}$ (i.e., the PoF is almost equal to 1) provides no additional information. 
On the contrary, the uncertainty reduction is significant on the right plot when updating $g_2$.
It is reasonable to conclude that, for this example,
only evaluating $g_2$ and not $g_1$ would reduce the optimization cost by avoiding a computationally intensive simulation ($g_1$) and by reducing the increase in data set size, while not impacting the algorithm performance.

It should be stated that the constraint selection idea that we just introduced in the context of robust optimization
would also stand for deterministic problems where the constraints are jointly modeled. 
In fact, it would be simpler: the one-step-ahead variance at $\mathbf x_{targ}$ of Equation \eqref{SC2_ConstSel} would be directly calculated without the integral on the $\mathbf u$'s. It would be minimized on $p$ only.

\section{Numerical implementation}
\label{sec:algosImplement}

\subsection{Algorithms compared}
\label{sec:theAlgos}
Four algorithms that correspond to the main combinations of the ideas presented above, i.e., coupled modeling of the constraints or not and constraint selection or not, will be compared:
\begin{itemize}
\item \textbf{REF}: separate model of each constraint, selection of a single joint value $(\mathbf x, \mathbf{u})$ at each iteration where the objective and all the constraints are evaluated. This is the reference method which was first proposed in \cite{elamri2021sampling}, hence the acronym.
\item  \textbf{SMCS}: Separate Model of the constraints and Constraint Selection. A single constraint and the associated $\mathbf u$ are identified at each iteration to provide the largest reduction in the probability of feasibility variance. 
\item \textbf{MMCU}: Multi-output Model of the constraints and Common value of $\mathbf u$. A single value of $\mathbf{u}$ is chosen at each iteration for which, after composition with $\mathbf x^{t+1}$, all the constraints are evaluated. 
\item \textbf{MMCS}: Multi-output Model and Constraint Selection. 
\end{itemize}
The REF algorithm could also have been called, with the keywords underlying the other acronyms, SMCU (Separate Model Common $\mathbf u$).
These four algorithms are described within a unified template in Algorithm~\ref{alg:template} and detailed in Appendix \ref{Alg}.
At the end of the four algorithms, the solution returned is the predicted solution to Problem \eqref{PBDef} as seen from the mean of the last GPs,
\begin{equation*}
\arg \min_{\mathbf{x}} m^{(t)}_Z(\mathbf{x}) ~\mbox{s.t.}~ \mathbb E[ C^{(t)}(\mathbf{x}) ] \leq 0
\end{equation*}
The last term, $\mathbb E[ C^{(t)}(\mathbf{x}) ] \leq 0$, is the expectation of satisfying the constraints with confidence $1-\alpha$.

\begin{algorithm}
\caption{Template of the robust constrained BO algorithms\\
{\small (the colors are associated to the variants: \textcolor{red}{coupled} or \textcolor{orange}{uncoupled} constraints, \textcolor{teal}{common $\mathbf u$} or \textcolor{green}{constraint selection})}}
\label{alg:template}
\begin{algorithmic}
\State Define initial data sets $\mathcal{D}_f^{(t_{init})},~ \mathcal{D}_{g_p}^{(t_{init})}, p=1,\ldots,l$ for the GPs of the objective function and constraints comprised of pairs $(\mathbf{x},\mathbf{u}) \in \mathcal{S_X} \times \mathcal{S_X}$ 
\State Evaluate the objective and the constraint functions through $f(\mathbf{x}, \mathbf{u}) ~,~ (\mathbf{x}, \mathbf{u}) \in \mathcal{D}_f^{(t_{init})}$ and $g_p(\mathbf{x}, \mathbf{u}) ~,~ (\mathbf{x}, \mathbf{u}) \in \mathcal{D}_{g_p}^{(t_{init})} ,~ p=1,\ldots,l$ $\Rightarrow \mathbf{f}^{(t_{init})}, \mathbf{g_1}^{(t_{init})}, \ldots , \mathbf{g_l}^{(t_{init})}$.
\State Let $t \gets  t_{init}$
\While{$t \leq t_{max}$}
\State Train the GPs, $F$ w.r.t. $(\mathcal{D}_f^{(t)},\mathbf{f}^{(t)})$ \\
\hspace{1cm} and \textcolor{red}{$\mathbf{G}$} w.r.t. $\left( (\mathcal{D}_{g_1}^{(t)} , \mathbf{g_1}^{(t)}), \ldots , (\mathcal{D}_{g_l}^{(t)} , \mathbf{g_l}^{(t)}) \right)$ \\
\hspace{1cm} in a \textcolor{red}{coupled} (\textcolor{red}{MM}\textcolor{teal}{CU}, \textcolor{red}{MM}\textcolor{green}{CS}) or \textcolor{orange}{independent} (\textcolor{orange}{R}\textcolor{teal}{EF}, \textcolor{orange}{SM}\textcolor{green}{CS}) way
\State Compute the targeted variable by maximizing the EFI (Equation \eqref{eq:EFI3}),\\
\hspace{1cm} $\mathbf{x}_{targ} = \mbox{arg} \min_{\mathbf{x}} \mathbb E[FI^{(t)}(\mathbf{x})]$
\State Set $\mathbf{x}^{t+1} = \mathbf{x}_{targ}$
\State Compute the random parameters :\\
\begin{itemize}
\item \textbf{either} (\textcolor{orange}{R}\textcolor{teal}{EF}, \textcolor{red}{MM}\textcolor{teal}{CU}) in a \textcolor{teal}{unified} fashion for the objective and constraints by minimizing the proxy to the one-step-ahead EFI variance (Equation \eqref{eq:proxy}), $\textcolor{teal}{\mathbf{u}^{t+1}} = \arg \min_\mathbf{u} S(\mathbf{x}_{targ},\mathbf{u}) $
\item \textbf{or} (\textcolor{red}{MM}\textcolor{green}{CS}, \textcolor{orange}{SM}\textcolor{green}{CS}) \textcolor{green}{separately} for the objective function and the selected constraint:
\begin{itemize}
\item \textcolor{green}{$\mathbf{u}_f^{t+1}$} minimizes the one-step-ahead improvement variance (Equation \eqref{eq:S1}), $\mathbf{u}_f^{t+1} = \mbox{arg} \min_{\mathbf{u}} S_f(\mathbf x_{targ},\mathbf u)$ 
\item $\mathbf u_g^{t+1}$ and the constraint index minimize the one-step-ahead variance of the feasibility probability (Equation \eqref{SC2_ConstSel}) , $\textcolor{green}{\{\mathbf u_g^{t+1}, p \}} = \mbox{arg} \min_{\mathbf{u},p'} S_g(\mathbf x_{targ},\mathbf u,p')$
\end{itemize}
\end{itemize}
\State Update $\mathcal{D}_f^{(t)},~ \mathcal{D}_{g_p}^{(t)}, p=1,\ldots,l$
\State Calculate $f(,)$ and (one -- \textcolor{red}{MM}\textcolor{green}{CS}, \textcolor{orange}{SM}\textcolor{green}{CS} -- or many -- \textcolor{orange}{R}\textcolor{teal}{EF}, \textcolor{red}{MM}\textcolor{teal}{CU} --)\\
\hspace{1cm} constraints at the new joined point(s), update $\mathbf{f}^{(t)}$ and $\mathbf{g_1}^{(t)},\ldots,\mathbf{g_l}^{(t)}$
\State $t \gets t+1$
\EndWhile\\
\Return the best feasible point according to the means of the processes,\\
\hspace{1cm} $\arg \min_{\mathbf{x}} m^{(t)}_Z(\mathbf{x}) ~\mbox{s.t.}~ \mathbb E[ C^{(t)}(\mathbf{x}) ] \leq 0$ 
\end{algorithmic}
\end{algorithm}

\subsection{Remarks about the implementation and code availability}
\label{ImplemNotes}
The algorithms presented and compared in this article are implemented in \texttt{R}. 
The scalar, independent, GPs are programmed thanks to the \texttt{DiceKriging} \cite{roustant2012dicekriging} toolbox. 
The multi-ouput GPs rely on the \texttt{kergp} \cite{deville2015package} toolbox which facilitates the implementation of the kernel seen in Equation \eqref{eq:kernelTensorP}. Please note that the source code allowing to reproduce the results presented in this paper on a simple 2-dimensional test case can be found on the following repository : \url{https://github.com/JPelamatti/Robust_Bayesian_Optimization}.\\
We now comment on their numerical implementations.

\subsubsection{Multi-output GP likelihood optimization}
When dealing with the multi-output model of the constraints, the conditioning data set is often made of the various constraints evaluated at the same location. 
That is, the data set may contain several points with the same $\mathbf{x}$ and $\mathbf{u}$ which only differ in $p$, the constraint indicator. 
This setting affects the algorithms which are eligible to maximize the GP likelihood. 
Indeed, the numerical computation of the likelihood gradient will return null gradient terms with respect to several GP hyper-parameters, namely the parameterization of the hypersphere kernel (Appendix \ref{sec:hypersphereKernel}) for the discrete indices $p$ contains $\theta_{i,j}$ hyper-parameters whose gradients vanish.
For this reason, the training of the multi-output GP models was performed with the COBYLA \cite{powell1994direct} optimization algorithm, which does not require gradients while still being quite efficient.
COBYLA searches were restarted 20 times from random initial points to make the likelihood maximization more robust.

\subsubsection{Other internal optimizations}
On top of the Gaussian process hyperparameter learning discussed above, 
the algorithms of this article include several internal optimizations such as 
the maximization of the feasible improvement (Equation \eqref{eq:EFI3}), 
the minimization of the one-step-ahead variance (Equation \eqref{eq:uplus1}) 
and the computation of a feasible mean performance (Equation \eqref{CurrFeasMin}).

Unlike standard Bayesian optimization, these auxiliary optimizations have a non-negligible computational cost because they require Monte Carlo simulations to estimate expected values and probabilities of the candidate solutions (this is for instance illustrated in Appendix \ref{sec:VarItp1calc}).

In order to afford a fair comparison between the algorithms studied, the results of Section \ref{NumRes} are obtained by optimizing the acquisition functions (Equations \eqref{eq:EFI3} and \eqref{eq:uplus1}) with respect to a predetermined set of candidate solutions randomly generated through a Latin Hypercube Sampling (LHS) procedure. The size of the considered candidate solutions set varies as $\simeq 500\times$dimension.
This implementation is not viable in larger dimensional problems (typically in dimension 5 and above), 
as the number of candidate solutions required to sufficiently fill the search space would be too large. 
In these cases, proper optimization algorithms (\emph{e.g.,} gradient-based, COBYLA, heuristic algorithms) are required to keep the number of candidate solutions evaluations at a computable level.

Notice also that the evaluations of the acquisition functions at different locations of the joint space are independent from each other. 
This property could be exploited to evaluate several candidate solutions in parallel, thus allowing toreduce the computation time required by the optimization of these acquisition functions.

\section{Experiments and results}
\label{NumRes}
This Section reports the experiments done with the four algorithms (REF, SMCS, MMCU, MMCS) on two analytical test cases as well as on a more realistic engineering design problem. 

\subsection{Experiments set-ups and performance metrics}
In order to take into account the variability coming from the initial data set, the runs are repeated 10 times with different initial data sets obtained through an optimized LHS, while the random variables are independently sampled from their specific distributions. 
The initial data set of the multi-constraints model is the union of the data sets for the scalar outputs. It is therefore $l$ times larger than each training set of the scalar GPs.

The same computational budget is allocated to every method. 
For comparing methods with and without constraint selection, the budget is expressed in terms of the number of constraints evaluations and an iteration is defined as $l$ evaluations of one constraint, regardless of repeated evaluations. 

At a given iteration, the performance of the algorithms is assessed in terms of the best current objective function value that is feasible with a probability of $1-\alpha$ where $\alpha=0.05$ or $\alpha=0.1$ depending on the considered problem.
This value is computed with respect to the true problem functions. 

On top of the previously discussed likelihood optimization, the BO algorithms studied embed auxiliary optimizations: the identification of the current feasible minimum, the optimization of the expected feasible improvement and the minimization of the one-step-ahead variance.
To limit the impact of the randomness in these auxiliary calculations on the tests, the auxiliary optimization are performed by determining the optimum among a fixed set of points. A different set of points is used for each optimization repetition.

\subsection{2-dimensional analytical test case}
The first test case depends on 1 design variable and 1 random variable, and is subject to two constraints: 
\begin{equation}
\begin{split}
&\min_{\mathbf{x}} \mathbb E_U[f(\mathbf{x},\mathbf{U})] \ \ s.t. \ \ \mathbb P(g_p(\mathbf{x},\mathbf{U}) \leq 0, p = 1,2) \geq 0.95 ~, \\
&\text{where }\quad   f(x,u) = (x-10)^3 + (u-20)^3 ~,\\ 
& \phantom{\text{where }\quad} g_1(x,u)  = -(x-5)^2 - (u-5)^2 + 500~, \\ 
& \phantom{\text{where }\quad} g_2(x,u)  = (x-6)^2 + (u-5)^2 - 9000~, 
\end{split}
\label{eq:2danalytic}
\end{equation}
$\mathbf{x} \in [13,100]$ and $\mathbf{U} \sim \mathcal{U}([0,100])$. 
This test case is intended to highlight the potential advantage of multi-output GP modeling thanks to the large correlation between the constraints. 
Additionally, given the expression of the two constraints, the covariance between $p=1$ and $p=2$ within the multi-output GP is expected to be negative. 
The discrete kernel parameterization should therefore be able to return negative values, which happens with the hypersphere decomposition.

\begin{figure}[h!]
\includegraphics[width=\textwidth]{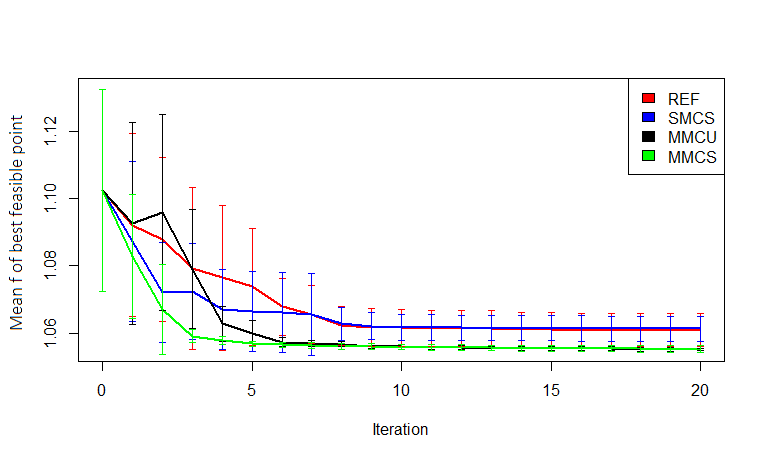}
\caption{Convergence of $\mathbb E_{\mathbf U} f$ for the best feasible point of the 4 optimization algorithms on the 2-dimensional Problem \eqref{eq:2danalytic}.
Median and 25\% and 75\% quartiles estimated from 10 repetitions.}
\label{Conv2d}
\end{figure}

\begin{figure}[h!]
\includegraphics[width=\textwidth]{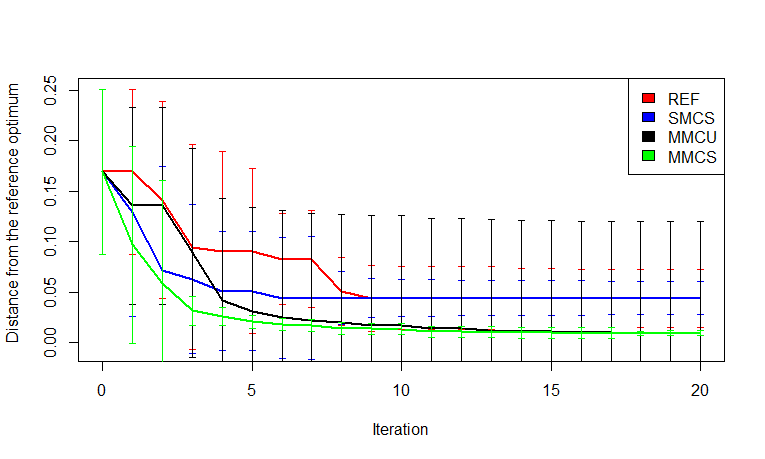}
\caption{Distance between the current best point and the optimum as a function of the iteration, for the 4 compared BO algorithms on the 2-dimensional Problem \eqref{eq:2danalytic}. Median and 25\% and 75\% quartiles estimated from 10 repetitions.
}
\label{ResDistance2d}
\end{figure}

\begin{figure}[h!]
\includegraphics[width=\textwidth]{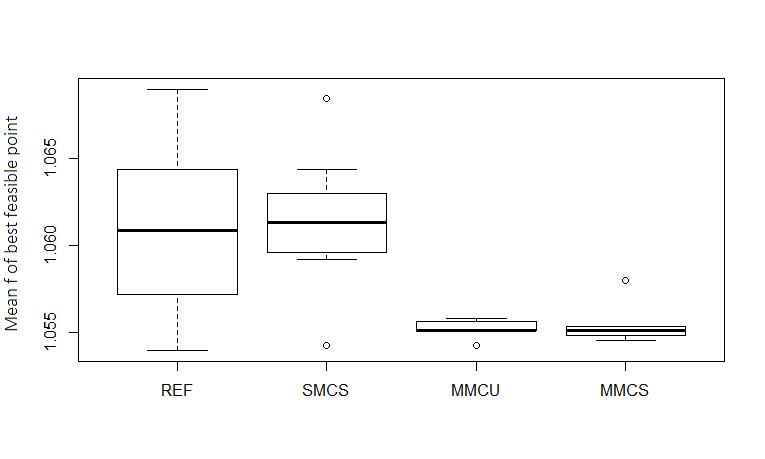}
\caption{Summary statistics for the best mean objective, $\mathbb E_{\mathbf U} f$, of feasible points obtained by the 4 optimizers on the 2-dimensional Problem \eqref{eq:2danalytic}, 10 repetitions of the runs.
}
\label{Res2d}
\end{figure}

The GP models of the constraints and the objective function are initialized with a data set of 6 samples per function. 
A budget of 40 constraint evaluations is allocated to the optimization iterations. 
The convergence rates of the 4 algorithms compared over the 10 repetitions are shown in Figure \ref{Conv2d}. An increase in the median best feasible value is seen around the third iteration for the MMCU algorithm. This can happen when some of the runs do not find a feasible solution at the beginning and are thus not taken into account. When an additional feasible point arises, the median best feasible objective may increase.
A complementary view on the convergence expressed in terms of the distance between 
the current best location and the true optimum is provided in Figure \ref{ResDistance2d}. 
The reference (true) optimum was found by exhaustive enumeration on the true functions and constraints.
The best mean objectives of feasible solutions found by the algorithms over the 10 repetitions are summarized in Figure \ref{Res2d}.

\begin{figure}[h!]
\includegraphics[width=\textwidth]{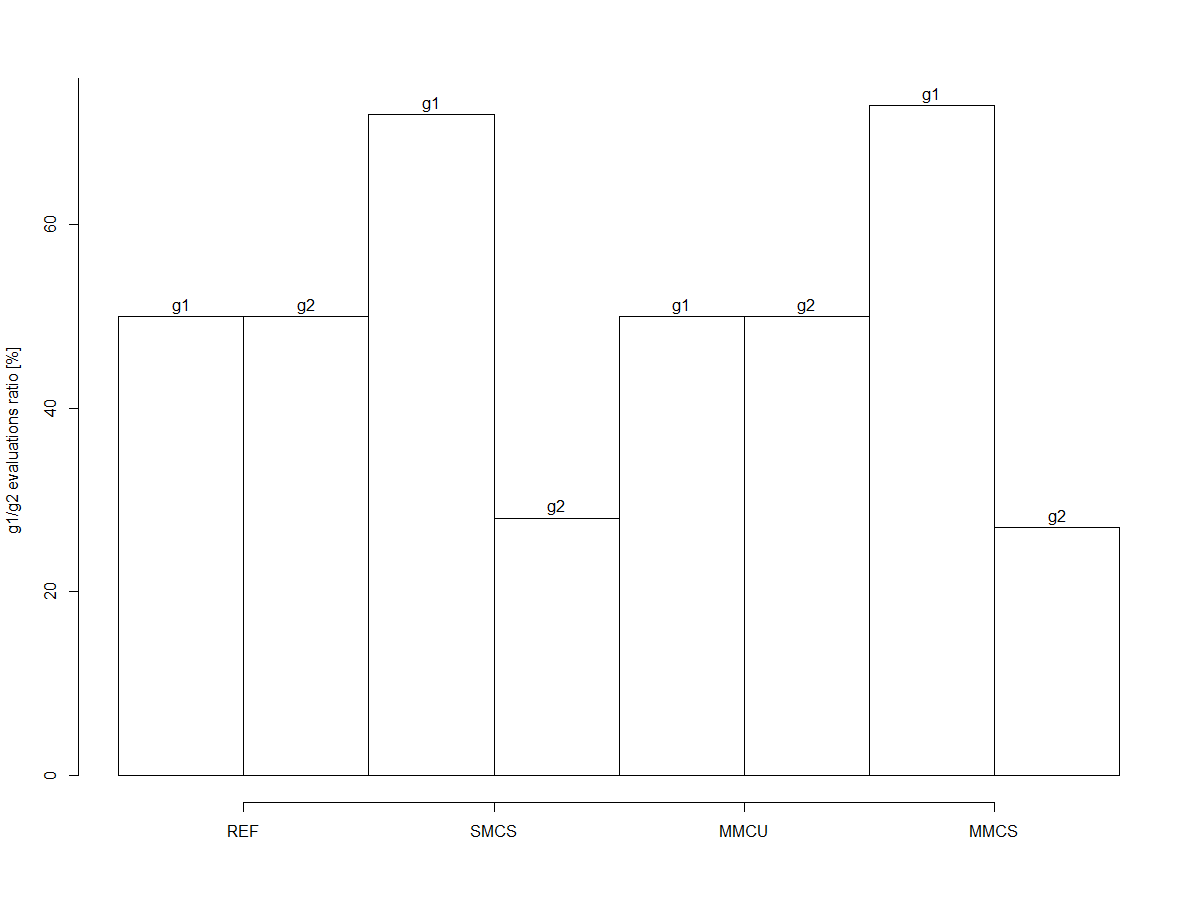}
\caption{Proportions of $g_1(\cdot)$ and $g_2(\cdot)$ evaluations made by the 4 algorithms on the 2-dimensional Problem \eqref{eq:2danalytic}, averaged over 10 repetitions.
}
\label{Constratio2d}
\end{figure}

\begin{figure}[h!]
\includegraphics[width=\textwidth]{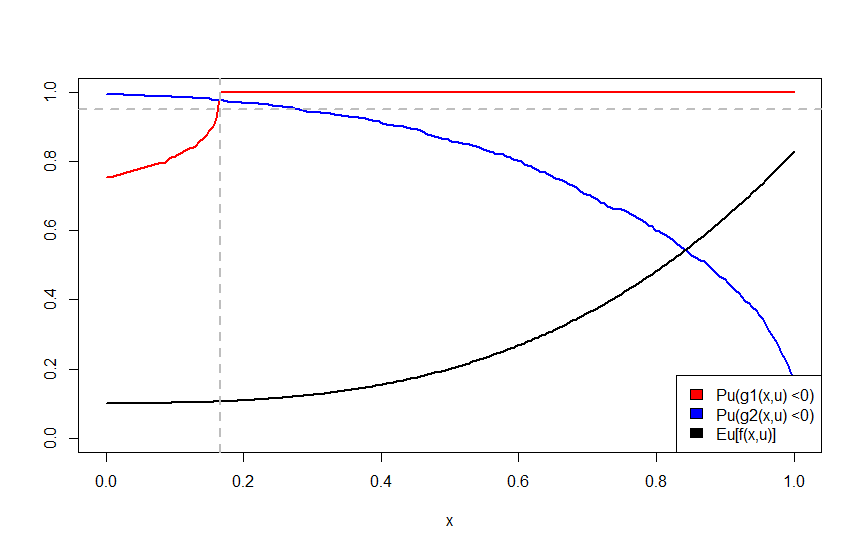}
\caption{Profiles of $P_U(g_1(\mathbf{x},\mathbf{u}) \leq 0)$, $P_U(g_2(\mathbf{x},\mathbf{u}) \leq 0)$ and $E_U(f(\mathbf{x},\mathbf{u}))$. 
The optimum $\mathbf x^\star$ is indicated by the vertical grey dashed line.
$E_U(f(\mathbf{x},\mathbf{u}))$ is rescaled for readability.
}
\label{PoF2d}
\end{figure}

For this test case, the 3 methods introduced here, SMCS, MMCU and MMCS, converge faster during the first iterations than the reference method (REF).  It is also observed that constraint selection, as exercised by SMCS and MMCS, accelerates the initial convergence (MMCS is better than MMCU).
The effect of constraint selection is detailed in Figure \ref{Constratio2d} which shows an histogram of the frequency with which each constraint is called. 
By construction, both the REF and MMCU approaches evaluate $g_1(\cdot)$ and $g_2(\cdot)$ the same number of times, i.e., they both represent 50\% of the overall constraints calls. 
On the contrary, both the SMCS and MMCS approaches evaluate $g_1(\cdot)$ more often than $g_2(\cdot)$, 70\% against about 30\% of the time.
An explanation for why the infill procedure prioritizes $g_1(\cdot)$ over $g_2(\cdot)$ is found in Figure \ref{PoF2d} 
which shows the profiles of $\mathbb P_U(g_1(\mathbf{x},\mathbf{u}) \leq 0)$, $\mathbb P_U(g_2(\mathbf{x},\mathbf{u}) \leq 0)$ 
and $\mathbb E_U(f(\mathbf{x},\mathbf{u}) \leq 0)$. 
It is seen that at the constrained optimum $\mathbf x^\star$ (marked by the vertical grey dashed line), $g_1(\cdot)$ is saturated 
while $g_2(\cdot)$ is above the threshold $1-\alpha$ (marked by the horizontal grey dashed line). 
As a consequence, the surrogate of $g_1(\cdot)$ is required to be more precise than that of $g_2(\cdot)$ to accurately assess feasibility
as measured, in our approaches, by the variance of the feasibility probability through $S_g$.

Later in the run, as seen in Figures \ref{Conv2d} and \ref{ResDistance2d}, the multi-output model enables to better refine the best solution, thus MMCS and MMCU consistently converge to feasible points with lower mean objective than the algorithms with independent constraint models, REF and SMCS.

\begin{figure}[h!]
\includegraphics[width=\textwidth]{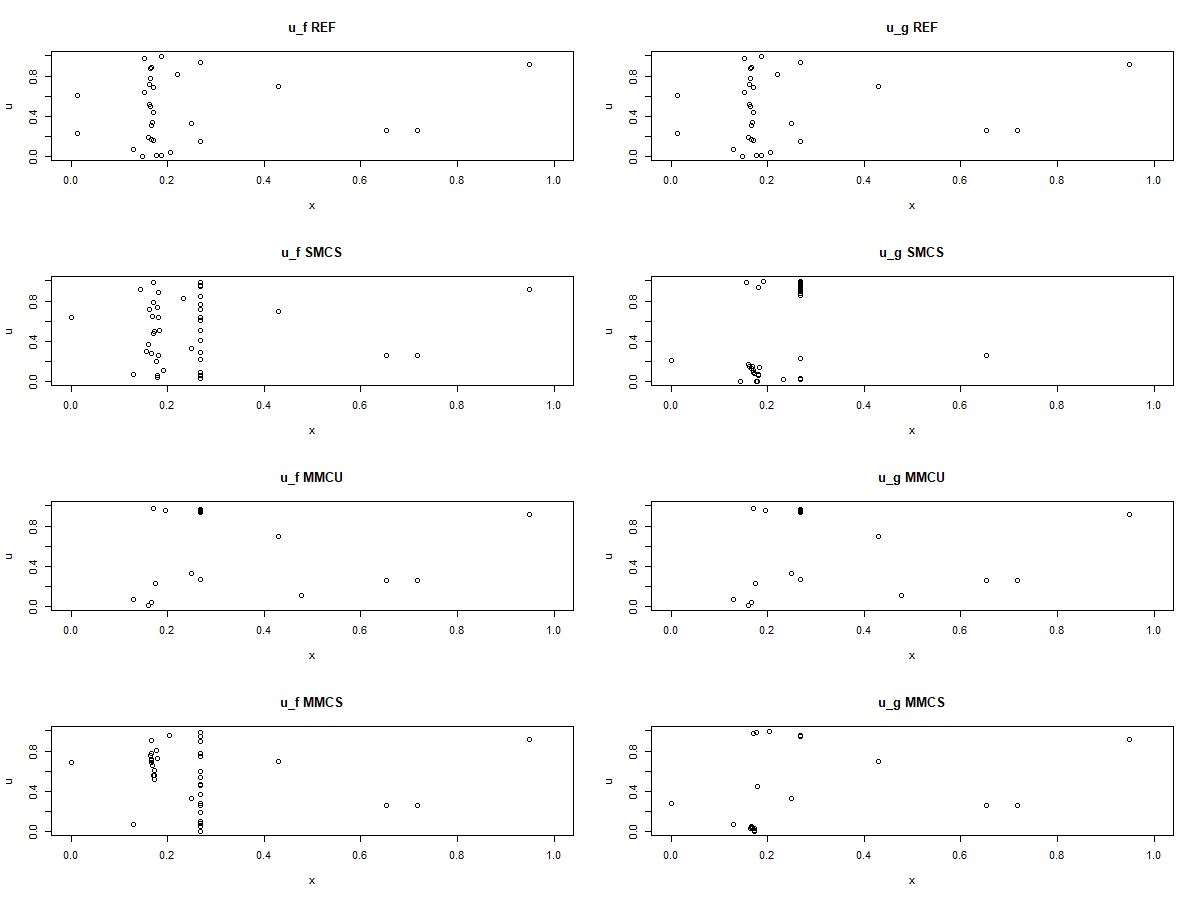}
\caption{
DoE's generated by one execution of the different algorithms on the 2-dimensional Problem \eqref{eq:2danalytic}.
$(x^i,u^i)~,~i=1,\ldots,40$ for the objective function (left column) and the constraints (right column), for REF, SMCS, MMCU and MMCS from row 1 to 4, respectively. }
\label{Infilled2d}
\end{figure}

The benefit of selecting independent $\mathbf{u}$'s for the evaluation of the objective and the constraint functions 
is illustrated in Figure \ref{Infilled2d} where the DoEs generated by one execution of the different algorithms is plotted.
It is observed that the algorithms that sample independently $u$ for the objective and the constraint functions (SMCS and MMCS) have the objective function samples $u_f$ that explores the entire search space while the constraints samples $u_g$ cluster on the boundaries. 
Instead, by construction, $u_f=u_g$ for REF and MMCU.
For all algorithms, the distribution of the variables $x$ has 2 clusters, one in the optimal region, $x^\star \approx 0.17$, 
and one where the second, inactive constraint, $0.95-\left(\mathbb P_U(g_2(x,U)<0\right)$, crosses 0 which happens for $x \approx 0.28$ (cf. Figure \ref{PoF2d}).

\subsection{4-dimensional analytical test case}
Our second test case has 2 design variables, 2 random variables, i.e., the models are built in 4 dimensions, and 2 constraints.
It was already introduced as Problem \eqref{eq:4d_analytic}.
It is also a modified version of a benchmark found in \cite{elamri2021sampling}, with the addition of a second constraint.

\begin{figure}[h!]
\includegraphics[width=\textwidth]{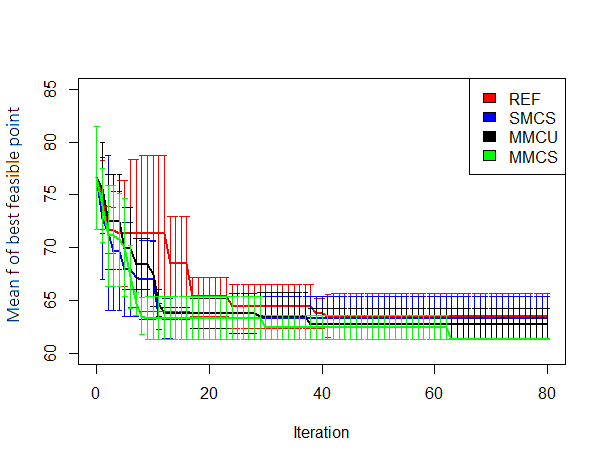}
\caption{Convergence rate (median and interquartile interval estimated from 10 runs) of the 4 BO algorithms on the 4-dimensional analytical test case \eqref{eq:4d_analytic}.
}
\label{Conv4d}
\end{figure}

\begin{figure}[h!]
\centering
\includegraphics[width=.7\textwidth]{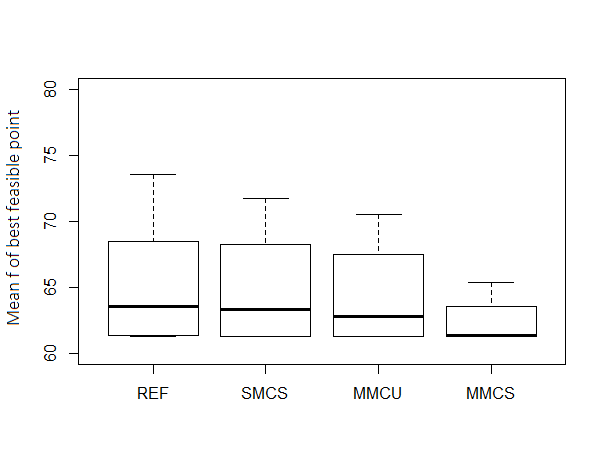}
\caption{Summary statistics over 10 repetitions of the objective, $\mathbb E_{\mathbf U}f(\widehat{\mathbf x^\star},\mathbf U)$ at the best feasible point $\widehat{\mathbf x^\star}$ obtained by the 4 BO algorithms on the 4-dimensional analytical test case \eqref{eq:4d_analytic}.}
\label{Final4d}
\end{figure}

The GP models of either the constraints or the objective functions are initialized with a data set of 30 points. 
Subsequently, a budget of 160 constraint evaluations is allocated to the optimization which, divided by 2 since there are 2 constraints, make 80 iterations.
The convergence rate of the 4 algorithms averaged over the 10 repetitions are shown in Figure \ref{Conv4d}. 
Summary statistics of the objective, $\mathbb E_{\mathbf U}f(\widehat{\mathbf x^\star},\mathbf U)$ at the best feasible point $\widehat{\mathbf x^\star}$ found by the algorithms over the 10 repetitions are provided in Figure \ref{Final4d}.
These results show that, similarly to the 2-dimensional test case, the three proposed variants (SMCS, MMCU, MMCS) 
outperform the reference algorithm (REF), both in terms of convergence speed and final performance of the solution while the interquartile distances remain inferior to that of REF.
\begin{figure}[h!]
\includegraphics[width=\textwidth]{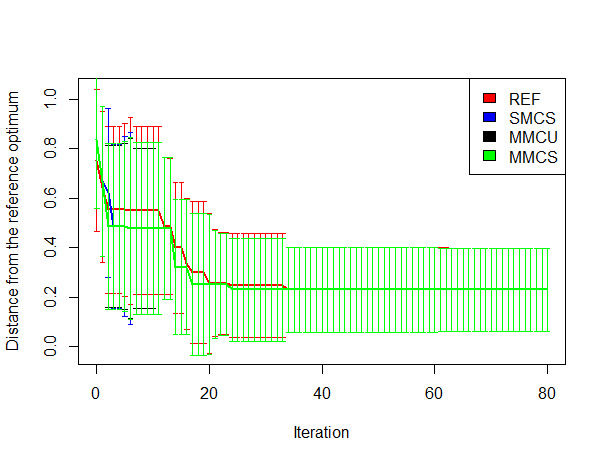}
\caption{Distance (median and interquartile interval) between the best solution and the reference solution obtained by the 4 BO algorithms on the 4-d analytical test case \eqref{eq:4d_analytic}.}
\label{ResDistance4d}
\end{figure}
The convergence rate of the algorithms is expressed in terms of Euclidean distance between the current best point and the reference optimum
in Figure \ref{ResDistance4d}.
This reference optimum was found by exhaustive enumeration on the true functions and constraints.
With this metric, the difference in performance between the various methods is less noticeable, but the proposed variants still show a faster convergence rate compared to REF. 
The MMCS variant provides the best results according to both metrics of comparison.

\begin{figure}
\centering
\begin{subfigure}{0.6\textwidth}
\includegraphics[width=\textwidth]{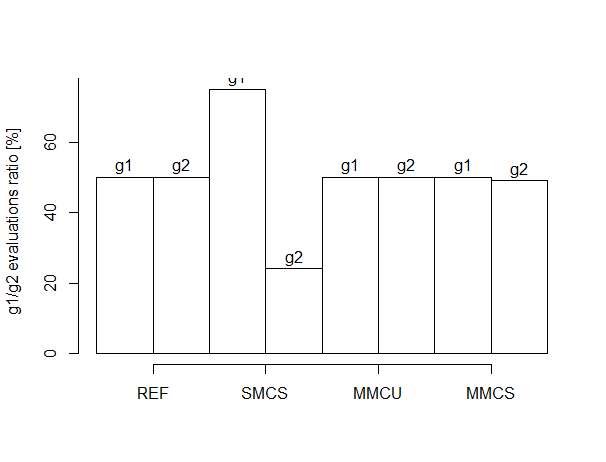}
\caption{Proportions of $g_1(\cdot)$ and $g_2(\cdot)$ evaluations averaged over 10 repetitions.}
\label{Const4d}
\end{subfigure}
\\
\begin{subfigure}{0.6\textwidth}
\includegraphics[width=\textwidth]{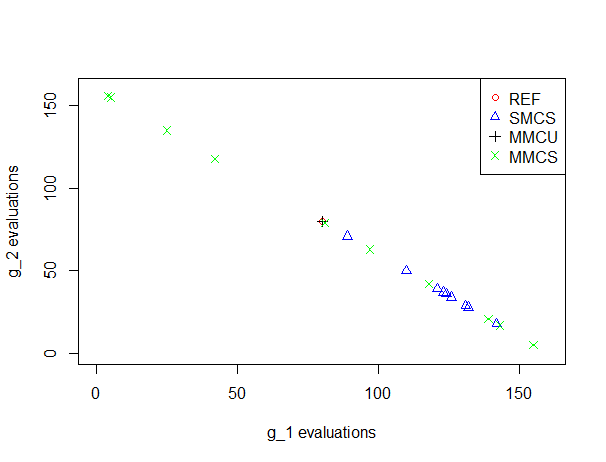}
\caption{Details of the number of calls to $g_1(\cdot)$ and $g_2(\cdot)$ for each of the 10 repetitions.}
\label{ConstRep4d}
\end{subfigure}
\caption{Number of calls to $g_1(\cdot)$ vs. $g_2(\cdot)$ for the REF, SMCS, MMCU and MMCS algorithms, 4-dimensional Problem \eqref{eq:4d_analytic}.}
\label{fig:constraint_4d}
\end{figure}

In order to further analyze the difference between the methods, Figure \ref{fig:constraint_4d} shows how the evaluations are distributed between the two constraints.
Whereas the SMCS variant evaluates more often, on the average of the 10 runs, $g_1$ than $g_2$, the MMCS variant seems to evenly calculate the two constraints. 
The number of evaluations is detailed for each run in Figure \ref{ConstRep4d}. 
It can be seen that the apparent equal evaluations of MMCS are due to an averaging effect: the MMCS focuses its computational effort on either one of the two constraints depending on the run. 
This is because the two constraints are strongly correlated and both are active. Therefore, either one of the two constraints provides a useful information from an optimization perspective which can be modeled by the multi-output GPs and exploited by the constraint selection mechanism.

\subsection{Industrial test case: compressor rotor design}
\begin{figure}[h!]
\centering
\includegraphics[width=0.8\textwidth]{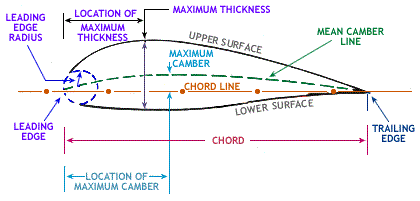}
\caption{Schematic representation of the rotor test case and its design variables, courtesy of Safran Tech.}
\label{Schema_Rotor}
\end{figure}

Our last test case is the design of an airplane engine compressor rotor. 
More specifically, the NASA rotor 37, an axial flow compressor originally designed and experimentally tested by Reid and Moore in 1978 \cite{reid1978design,dunham1998cfd} is considered. 
The objective of the design is to maximize the compressor polytropic efficiency. 
The design concerns the shape of the rotor, as described by 5 quantities sketched in Figure \ref{Schema_Rotor}:
the chord, the maximum thickness, the location of the maximum thickness, the pitch and the sweep angles. 
Each shape feature evolves along the rotor axis and is parameterized with the help of a spline function with 4 control points, thus resulting in 20 actual design variables.
Additionally, the rotor performance depends on 7 random parameters $\mathbf U$:
the tip gap and the blade rugosity (manufacturing uncertainties), the inflow pressure, the temperature 
and the azimuthal momentum (inflow uncertainties), the flow rate and the rotation speed (operational uncertainties).  
Finally, the optimization problem is subject to 5 constraints related to the inlet and outlet relative flow angles ($g_1,g_2$), the deceleration of the flow to prevent separation ($g_3$), the Lieblin coefficient to avoid excessive loading  ($g_4$) and the Mach number of the rotor tip  ($g_5$).

The rotor design problem is formulated as follows:
\begin{equation}
\min_{\mathbf{x}} \mathbb E_U[f(\mathbf{x},\mathbf{U})] \ \ s.t. \ \ \mathbb P(g_p(\mathbf{x},\mathbf{U}) \leq 0, p = 1,\ldots, 5) \geq 0.9
\end{equation}
with $\mathbf{x} \in \mathcal{S_X} \subset \mathbb{R}^{20}$ and $\mathbf{U} \sim \mathcal{U}([\text{LB},\text{UB}]^7)$, 
where $\text{LB}$ and $\text{UB}$ represent the lower and upper bounds of the uniformly distributed random parameters.

With 27 variables and parameters, and 5 constraints, the joint space of the rotor design problem is considerably larger than those of the previous analytical test cases. 
Robust optimizations are computationally much more expensive.
As a result, various adjustments in the experiments have been deemed necessary.
First, the reliability in constraint satisfaction has been set to 90\% ($\alpha = 0.1$).
Second, only a single compressor optimization is performed, i.e., the tests are not repeated. 
The results provided hereafter should be taken with some caution as the variability caused by the initial DoE and the auxiliary optimizations (GP learning, acquisition criteria optimizations) is not studied.
Third, only two of the four methods are considered: REF and MMCS. 

\begin{figure}[h!]
\includegraphics[width=\textwidth]{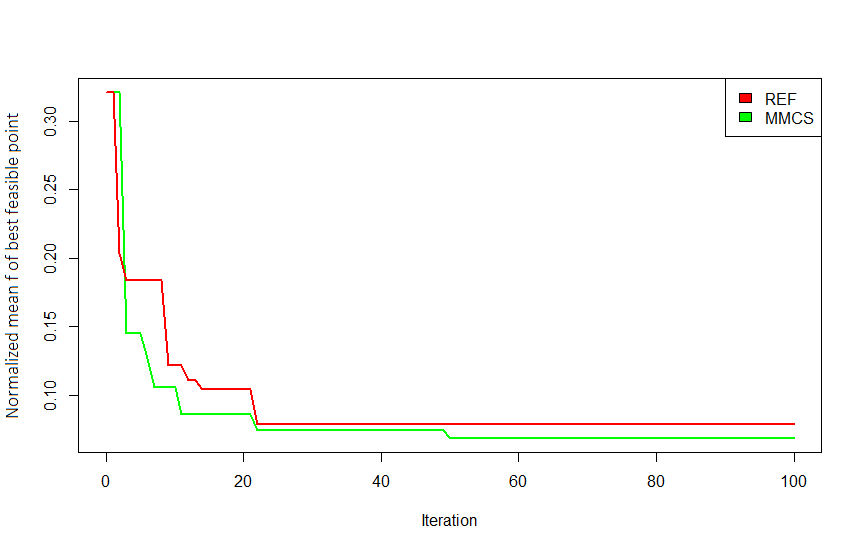}
\caption{Convergence of the best feasible $\mathbb E_{\mathbf U} f$ for the REF and MMCS algorithms, compressor rotor design test case (1 execution).
}
\label{ConvSafran}
\end{figure}

The GP models are trained with initial DoEs of 100 data samples per function. 
With 5 constraints, this results in the initial multi-output GP of MMCS being conditioned on 500 data samples. 
Subsequently, 500 constraint functions evaluations are performed for the actual optimization (i.e., 100 iterations). 
The convergences of the REF and MMCS algorithms are plotted in Figure \ref{ConvSafran}.
This result shows that the MMCS variant, in green, with its multi-output constraints model and constraint selection, 
outperforms the reference (REF) algorithm, in red, as it provides a faster initial convergence speed 
as well as a better final mean objective. 
In order to further investigate the effects of selecting the single most informative constraint at every iteration, 
the number of evaluations of each constraint is recorded in Table \ref{TableCompEval}.
\begin{table}[h!]
\centering
\begin{tabular}{|l|c|c|c|c|c|}
\hline
Constraint & $g_1$ & $g_2$ & $g_3$ & $g_4$ & $g_5$ \\  \hline
N$^\circ$ evaluations REF & 100 & 100 & 100 & 100 & 100  \\  \hline
N$^\circ$ evaluations MMCS & 103 & 202 & 17 & 98 & 80 \\  \hline 
\end{tabular}
\caption{Number of evaluations of each constraint in the compressor rotor design test case carried out by the REF and MMCS methods.}
\label{TableCompEval}
\end{table}
It is seen in the Table that the MMCS focuses a larger part of the computational budget on some constraints, in particular $g_2$ and to a smaller extent $g_1$, and tends to neglect other constraints (e.g., $g_3$).

\section{Conclusion and perspectives}
\label{Conclusions}
In this paper, three compatible extensions to the robust Bayesian optimization algorithm were proposed in order to solve computationally intensive chance constrained problems. 
These extensions are based 
\textit{i)} on the creation of a coupled model of the constraints,
\textit{ii)} on the independent selection of the constraint and the associated random parameter where it will be evaluated. This was made possible thanks to a decomposition of the acquisition criterion and to an output-as-input scheme.
The results obtained on both academic and an engineering design test cases show 
that the two changes contribute to a convergence which is both more accurate and more stable with respect to the initial data set when compared to the original method. 
We therefore think that coupled constraint models and constraint selection are key features of BO methods 
when dealing with costly design problems affected by the presence of random variables.

Additional testing might be necessary to better assess the effects and possible limitations of the proposed approaches. 
For example, it could be interesting to assess whether relying on multi-output models can induce a loss of performance when the constraints are not correlated. 

From a methodological perspective  the MM** algorithms should be further developed to allow the identification and the processing of groups of correlated constraints.
Finally, another continuation to this work would be to model the correlation between the objective and the constraint functions, 
with the anticipated challenge of having to define an acquisition criterion fitted to this assumption as well as computationally tractable.

\subsubsection*{Acknowledgements:} This work was partly supported by the OQUAIDO research Chair in Applied Mathematics.

\newpage

\bibliographystyle{unsrt}
\bibliography{biblio}

\appendix

\section{Main notations and acronyms}
\label{sec:notations}
\begin{multicols}{2}
{\small (by alphabetical order)} 
\begin{itemize}[label={},itemsep=0mm]
\item BO~: Bayesian Optimization.
\item $c(.), C(.)$~: $1-\alpha -$ the probability of satisfying all the constraints, resp. its approximation as a process. Negative values mean feasibility with confidence $1-\alpha$.
\item $\mathcal D^{(t)},\mathcal D_f^{(t)},\mathcal D_{g_p}^{(t)}, \mathcal D_{\mathbf g}^{(t)}$~: Design of Experiment at iteration $t$, containing a collection of points $(\mathbf x^i,\mathbf u^i)$. The subscript $f$, $g_p$, $\mathbf g$ refer to the objective function, the $p$-th constraint or the joint constraints vector. 
\item DoE~: Design of Experiment.
\item I, FI, EI, EFI~: Improvement, Feasible Improvement, Expected Improvement, Expected Feasible Improvement.
\item $f(.,.), F(.,.)$~: function $\mathcal{S_X} \times \mathcal {S_U} \to \mathbb R$ and associated GP.
\item $\mathbf f^{(t)}$~: vector of evaluations of the function $f$ evaluated at the DoE $\mathcal D^{(t)}$ or $\mathcal D_f^{(t)}$
\item $g_p(,)$~: $p$-th constraint function $\mathcal{S_X} \times \mathcal {S_U} \to \mathbb R$, $p=1,\ldots,l$.
\item $\mathbf g^{(t)} , \mathbf g_p^{(t)} $~: vector of evaluations of the constraints functions $g_1(,) , \ldots , g_l,)$ (resp. only $g_p(,)$) evaluated at the DoE $\mathcal D^{(t)}$ or $\mathcal D_{\mathbf g}^{(t)}$ or $\mathcal D_{g_p}^{(t)}$.
\item $\mathbf g(.,.), \mathbf G(.,.)$~: vector of $l$ constraints and associated (not necessarily independent) GPs. 
\item $GP$~: Gaussian Process.
\item $k(),\mathbf K()$~: covariance function and matrix of a GP. Additional subscripts indicate which GP is concerned.
\item $l$~: number of constraints.
\item $m(),\mathbf m()$~: mean of a GP or of a vector of GPs. Additional subscript indicates which GP is concerned.
\item MMCS~: Multiple Model of the constraints and Constraint Selection algorithm.
\item MMCU~: Multiple Model of the constraints and Common $\mathbf u$ algorithm.
\item REF~: reference algorithm (independent constraints, same iterate for all constraints and objective function).
\item $S(.,.)$~: proxy to the one-step-ahead EFI variance.
\item $S_f(.,.)$~: one-step-ahead improvement variance.
\item $S_g(.,.)$~: one-step-ahead variance of the feasibility probability.
\item SMCS~: Single Models of the constraints and Constraint Selection algorithm.
\item $x,x',x^\star, x_{targ}$~: vectors of optimization variables $\in \mathcal{S_X} \subset \mathbb R^d$.
\item $u,u',\tilde u$~: vector of random parameters $\in \mathcal{S_U} \subset \mathbb R^m$.
\item $z(.), Z(.)$~: mean of $f(.,U)$ over $U$ and associated Gaussian process.
\end{itemize}
\end{multicols}

\section{Covariance kernel for a nominal input: the hypersphere decomposition.}
\label{sec:hypersphereKernel}
The hypersphere decomposition \cite{zhou2011simple} is a possible choice to parameterize a discrete covariance kernel. 
The underlying idea is to map each of the $l$ levels of the considered discrete variable onto a distinct point on the surface of a $l$-dimensional hypersphere:
\begin{eqnarray*}
&& \phi(z) : \mathcal S_z \rightarrow \mathbb{R}^l \\ 
&& \phi(z = z_m) = \sigma_z [ b_{m,0}, b_{m,1}, \dots , b_{m,l}]^\top \nonumber \quad \mbox{ for } m=1,\ldots , l
\end{eqnarray*}
where $b_{m,d}$ represents the $d$-th coordinate of the $m$-th discrete level mapping, and is calculated as follows:
\begin{align*}
& b_{1,1} = 1 & 
\\
& b_{m,d} = \cos \theta_{m,d} \prod_{k = 1}^{d-1} \sin \theta_{m,k} & \mbox{ for } d = 1,\dots,m-1 \\
& b_{m,m} =\prod_{k = 1}^{d-1} \sin \theta_{m,k} & \mbox{ for }  m \neq 1 \\
& b_{m,d} = 0  & \mbox{ for } d \geq m \neq 1
\end{align*}
with $-\pi \leq \theta_{m,d} \leq \pi$. It can be noticed that in the equations above, some of the mapping coordinates are arbitrarily set to 0. This allows to avoid rotation indeterminacies (\emph{i.e.,} an infinite number of hyperparameter sets characterizing the same covariance matrix), while also reducing the number of parameters required to define the mapping. The resulting kernel is then computed as the Euclidean scalar product between the hypersphere mappings, 
\begin{equation*}
k(z,z') = \phi(z)^\top\phi(z')~.
\end{equation*}
The discrete kernel can then be characterized as an $l\times l$ symmetric positive definite matrix $\mathbf{T}$ containing the covariance values between the discrete variable levels computed as:
\begin{equation*}
\mathbf{T} = \mathbf{L}^\top\mathbf{L}
\end{equation*}
where each element of $\mathbf{L}_{i,j}$ is computed as $b_{i,j}$,
\begin{equation*}
\mathbf{L} = \sigma_z \left[
\begin{array}{c c c c c}
1 & 0 & \dots & \dots & 0 \\
\cos \theta_{2,1}  & \sin \theta_{2,1} & 0 & \dots & \dots \\
\vdots & \vdots & \vdots & \vdots & \vdots \\
\cos \theta_{l,1} & \sin\theta_{l,1} \cos\theta_{l,2} & \dots & \cos \theta_{l,l-1} \prod_{d = 1}^{l-2} \sin\theta_{l,d}  & \prod_{d = 1}^{l-1} \sin\theta_{l,d}
\end{array} \right] ~.
\end{equation*}

In the output-as-input model of the constraints, Equation \eqref{eq:kernelTensorP}, the matrix $\mathbf T$ contains the covariance terms related to the constraint index, $k_p(i,j) = \mathbf T_{ij}~,~i,j = \{1,\ldots,l\}$.

\section{Probability of feasibility with dependent constraints}
\label{sec:detailsFeasImprov}
The probability of satisfying the coupled constraints intervenes in Equation \eqref{eq:EFI3} for the EFI acquisition criterion.
It can be estimated with GPs as:
\begin{equation}
\mathbb{P}(C^{(t)}(\mathbf{x}) \leq 0) \approx \frac{1}{N} \sum \limits_{k=1}^{N} \indic_{\big(1 - \alpha - \frac{1}{M} \sum \limits_{j=1}^{M}\indic_{\big(\mathbf{G}^{(t)}(\mathbf{x},\mathbf{u_j},\omega_k) \leq \mathbf{0} \big)} \leq 0\big)}.
\end{equation}
In practice, a set of $M$ instances of the random variables $\mathbf{U}$ is sampled from $\rho_U$. 
Subsequently, $N$ independent multi-output trajectories of $\mathbf{G}(\mathbf{x},\cdot)$ are simulated at the aforementioned sampled random variables. The probability of feasibility can finally be computed by simply counting the number of trajectories for which the ratio of samples associated to feasible constraints is larger than $1-\alpha$. 
Note that a multi-output GP prediction of the constraint vector is defined as feasible when all of its components are below their specific threshold value (which is $0$ in this work).

In \eqref{eq:EI}, the improvement is computed with respect to the  incumbent best feasible solution $z_{min}^{feas}$. 
However, the objective function mean is not observed
therefore $z_{min}^{feas}$ cannot  be  read in the data set used to condition the GPs $F$ and $\mathbf{G}$. For this reason, in \cite{elamri2021sampling} the incumbent best feasible solution is defined by taking into account the mean of the GP $Z^{(t)}(\mathbf{x})$ and the expected value of the process $C^{(t)}(\mathbf{x})$:
\begin{equation}
\label{CurrFeasMin}
z_{min}^{feas} =  \min_{\mathbf{x}} m^{(t)}_Z(\mathbf{x}) ~\mbox{s.t.}~ \mathbb E[ C^{(t)}(\mathbf{x}) ] \leq 0~.
\end{equation}
Given that the Fubini condition holds since the value of $C(\cdot)$ is bounded by definition, the expected value of the process $C$ can be written
\begin{eqnarray}
\mathbb E[C^{(t)}(\mathbf{x})] & =& 1 - \alpha - \mathbb E_\mathbf{U}[\mathbb E[1_{\mathbf{G} (\mathbf{x}, \mathbf{U}) \leq \mathbf{0} } ]] \nonumber\\
& =&  1 - \alpha - \int_{\mathbb{R}^m}
 \mathbb P(\mathbf{G} (\mathbf{x}, \mathbf{u}) \leq \mathbf{0}))\rho_\mathbf{U} (\mathbf{u}) d\mathbf{u}\nonumber \\
&= &  1 - \alpha - \int_{\mathbb{R}^m}
\Phi\left( \mathbf{0} -\mathbf{m}_\mathbf{G}^{(t)}(\mathbf{x},\mathbf{u}),\mathbf{K}_\mathbf{G}^{(t)}(\mathbf{x},\mathbf{u}) \right)\rho_\mathbf{U} (\mathbf{u}) d\mathbf{u} 
\end{eqnarray}
where $\Phi(\cdot)$ is the cumulative distribution function of a multivariate Gaussian distribution which, like the univariate version, is estimated numerically \cite{genz2009computation}. 
The above Equation and Equation \eqref{CurrFeasMin} are used to compute the current feasible minimum. 
If Equation \eqref{CurrFeasMin} yields no solution, $m^{(t)}_Z(\mathbf{x})$ at the  $\mathbf{x}$ providing the largest probability of feasibility is taken as the incumbent optimal solution.

\section{Estimating the proxy of the one-step-ahead feasible improvement variance}
\label{sec:VarItp1calc} 

During the optimization process, the value of the coordinates $\mathbf{u}^{t+1}$ of the point to be added to the training data set $\mathcal{D}$ is computed by minimizing the proxy of the one-step-ahead variance of the EFI at $\mathbf{x}_{targ}$ (see Equation \eqref{eq:uplus1}), which is defined as:

\begin{multline}\label{SC}
S(\mathbf{x}_{targ},\mathbf{\tilde{u}}) = Var\big( I^{(t+1)}(x_{targ})\big) \int_{\mathbb{R}^m} Var\big( 1_{\{\mathbf{G}^{(t+1)}(x_{targ},\mathbf{u}) \leq 0 \}} \big) \rho_{\mathbf{U}}(\mathbf{u}) d\mathbf{u},\\
 =~Var(\big(z_{\min}^{\mbox{feas}} - Z^{(t+1)}(x_{targ})\big)^+)\int_{\mathbb{R}^m} Var\big( 1_{\{\mathbf{G}^{(t+1)}(x_{targ},\mathbf{u}) \leq 0 \}} \big) \rho_{\mathbf{U}}(\mathbf{u}) d\mathbf{u}. \nonumber
\end{multline}

In this appendix, we recall some details about the computation of the first term in the previous equation. An expression of the improvement variance $Var\left( I^{(t+1)}(x_{targ})\right)$ which bears some similarities to the one of  expected improvement is given in \cite{elamri2021sampling} and can be expressed in terms of probability and density functions of Gaussian distributions:

\begin{eqnarray*}
EI^{(s)}(\mathbf{x})&=& \big(z_{\min}^\text{feas} - m_Z^{(s)}(\mathbf{x}) \big) \phi\bigg(\frac{z_{\min}^\text{feas} - m_Z^{(s)}(\mathbf{x})}{\sigma_Z^{(s)}(\mathbf{x})}\bigg) + \sigma_Z^{(s)}(\mathbf{x})  \Phi\bigg(\frac{z_{\min}^\text{feas} - m_Z^{(s)}(\mathbf{x})}{\sigma_Z^{(s)}(\mathbf{x})}\bigg),\\
&=&\psi_{EI}(m_Z^{(s)}(x), \sigma_Z^{(s)}(x)).
\label{firstterm}
\end{eqnarray*}

\begin{eqnarray*}
Var\left(I^{(s)}(\mathbf{x})\right)&=& EI^{(s)}(\mathbf{x}) \big(z_{\min}^\text{feas} - m_Z^{(s)}(\mathbf{x}) - EI^{(s)}(\mathbf{x})\big) \\
&&+ (\sigma_Z^{(s)})^2(\mathbf{x})  \Phi\bigg(\frac{z_{\min}^\text{feas} - m_Z^{(s)}(\mathbf{x})}{\sigma_Z^{(s)}(\mathbf{x})}\bigg),\\
&=&\psi_{VI}(m_Z^{(s)}(x), \sigma_Z^{(s)}(x)).
\end{eqnarray*}

At step $t+1$,  the training data set $\mathcal{D}^{(t)}$ is enriched by $(\tilde{\mathbf{x}}, \tilde{u})$ on which the  output $f(\tilde{\mathbf{x}},\mathbf{\tilde{u}})$ is unknown and represented by $F(\tilde{\mathbf{x}},\mathbf{\tilde{u}})$. 

As a consequence, $Var\left(I^{(t+1)}(x_{targ})\right)$  cannot directly be computed because of the randomness of  $m_Z^{(t+1)}(x_{targ})$. Indeed, 
$m_Z^{(t+1)}(x_{targ})$ is equal to 
 
\noindent $E\left(Z(x)| F(\mathcal{D}^{(t)})=f^{(t)}, F(\tilde{\mathbf{x}}, \tilde{u}) \right)$ and follows : 
\begin{equation*}
m_Z^{(t+1)}(x_{targ}) \sim \mathcal{N}\Bigg( m_Z^{(t)}(x_{targ}) , \bigg(\frac{\int_{\mathbb{R}^m} k_{*F}^{(t)}(x_{targ},\mathbf{u};\tilde{\mathbf{x}},\mathbf{\tilde{u}}) \rho_\mathbf{U}(\mathbf{u}) d\mathbf{u} }{\sqrt{k_{*F}^{(t)}(\tilde{\mathbf{x}},\mathbf{\tilde{u}};\tilde{\mathbf{x}},\mathbf{\tilde{u}})}}\bigg)^2 \Bigg). \\
\label{updatemeanZ}
\end{equation*}

We can note that  $\sigma_Z^{(t+1)}(x_{targ})$ is not random as it depends only on the location $(\tilde{\mathbf{x}},\mathbf{\tilde{u}})$ and not on the function evaluation at this point. By applying the law of total variance, it can be shown that:
\begin{eqnarray*}
Var \left( I^{t+1}(x_{targ})\right) & = & \mathbb E\left[Var\left(\big(z_{\min}^{\text{feas}} - Z(x_{targ})\big)^+ | m_Z^{(t+1)}(x_{targ}) \right)\right] \\
& &+Var\left[ \mathbb E\left(\big(z_{\min}^{\text{feas}} - Z(x_{targ})\big)^+ | m_Z^{(t+1)}(x_{targ}) \right)\right]. \\
&=& \mathbb E\left[ \psi_{VI}(m_Z^{(t+1)}(x), \sigma_Z^{(t+1)}(x)) \right] \\
& &+Var\left[ \psi_{EI}(m_Z^{(t+1)}(x), \sigma_Z^{(t+1)}(x)) \right]. \\
\label{firstterm2}
\end{eqnarray*}

This calculation  is performed numerically using samples of $m_Z^{(t+1)}(x_{targ}) $. For the sake of clarity, the reader is referred to the previous work (\cite{elamri2021sampling}) for details regarding the implementation.

\section{Flow charts of the algorithms}
\label{Alg}
This appendix contains the flowcharts of the four algorithms which are implemented and compared in the body of the article. 
The REF, SMCS, MMCU and MMCS methods are described in Algorithms \ref{alg:REF}, \ref{alg:SMCS}, \ref{alg:MMCU} and \ref{alg:MMCS}, respectively. 
REF stands for reference and was initially proposed in \cite{elamri2021sampling}. 
It is a method where the GPs of the constraints and objective function are independent, and the same pair $(\mathbf x^{t}, \mathbf u^{t})$ is added to every GP at each iteration. 
SMCS, which stands for Single Models of the constraints and Constraint Selection, has independent GPs, like the REF algorithm, but only one constraint is updated at each iteration. 
The random parameters of the objective function and the selected constraint, $\mathbf u_f$ and $\mathbf u_g$, are different.
MMCU means Multiple Model of the constraints and Common $\mathbf u$. The MMCU algorithm has a joint model of all the constraints and the same iterate $(\mathbf x^{t}, \mathbf u^{t})$ enriches all GPs. 
Finally, MMCS is the acronym of Multiple Model of the constraints and Constraint Selection. 
The MMCS algorithm relies on a joint model of the constraints and identifies at each iteration a single constraint and the associated random sample, $\mathbf u_g$, to carry out the next evaluation.

\begin{algorithm}
\caption{REF algorithm}
\label{alg:REF}
\begin{algorithmic}
\State Initialize $t_{init}$, $t_{max}$
\State Define the initial data set $\mathcal{D}^{(t_{init})} = \{ (\mathbf{x}^i,\mathbf{u}^i), \ i= 1,\ldots t_{init} \}$
\State Evaluate the objective and the constraint functions : 
\State $\mathbf{f}^{(t_{init})} = \left( f(\mathbf{x}^1, \mathbf{u}^1), \ldots , f(\mathbf{x}^{t_{init}}, \mathbf{u}^{t_{init}}) \right)$
\State $\mathbf{g}_1^{(t_{init})} = \left( g_1(\mathbf{x}^1, \mathbf{u}^1), \ldots , g_1(\mathbf{x}^{t_{init}}, \mathbf{u}^{t_{init}}) \right)$ 
\State \vdots
\State $\mathbf{g}_l^{(t_{init})} = \left( g_l (\mathbf{x}^1, \mathbf{u}^1), \ldots , g_l(\mathbf{x}^{t_{init}}, \mathbf{u}^{t_{init}}) \right)$ 
\State $t \gets t_{init}$
\While{$t \leq t_{max}$}
\State Train the GPs $F$ and $G_p$ for $p = 1, \ldots , l$ w.r.t. $\mathcal{D}^{(t)}$, $\mathbf{f}^{(t)}$ and $\mathbf{g}_p^{(t)}$
\State Compute the targeted variable by maximizing the EFI (Equation \eqref{eq:EFI3}),\\
\hspace{1cm} $\mathbf{x}_{targ} = \mbox{arg} \min_{\mathbf{x}} \mathbb E[FI^{(t)}(\mathbf{x})]$
\State Set $\mathbf{x}^{t+1} = \mathbf{x}_{targ}$
\State Compute the random parameter $\mathbf u^{(t+1)}$ as the minimizer \\
\hspace{1cm} of the proxy to the one-step-ahead EFI variance (Equation \eqref{eq:proxy}), \\
\hspace{1cm} $\mathbf{u}^{t+1} = \arg \min_\mathbf{u} S(\mathbf{x}_{targ},\mathbf{u}) $
\State $\mathcal{D}^{(t)} \gets \{\mathcal{D}^{(t)}, (\mathbf{x}^{t+1}, \mathbf{u}^{t+1}_f, \mathbf{u}^{t+1}_g)\} $
\State $\mathbf{f}^{(t)} \gets \{\mathbf{f}^{(t)}, f(\mathbf{x}^{t+1}, \mathbf{u}^{t+1})\} $
\State $\mathbf{g}_p^{(t)} \gets \{\mathbf{g}_p^{(t)}, g_p(\mathbf{x}^{t+1}, \mathbf{u}^{t+1})\} $ for $p = 1, \ldots, l$
\State $t \gets t+1$
\EndWhile \\
\Return the best feasible point according to the means of the processes,\\
\hspace{1cm} $\arg \min_{\mathbf{x}} m^{(t)}_Z(\mathbf{x}) ~\mbox{s.t.}~ \mathbb E[ C^{(t)}(\mathbf{x}) ] \leq 0$ 
\end{algorithmic}
\end{algorithm}

\begin{algorithm}
\caption{SMCS algorithm}
\label{alg:SMCS}
\begin{algorithmic}
\State Initialize $t_{init}$, $t_{max}$
\State Define the initial data set $\mathcal{D}^{(t_{init})} = \{ (\mathbf{x}^i,\mathbf{u}^i), \ i= 1,\ldots t_{init} \}$
\State Evaluate the objective and the constraint functions : 
\State $\mathbf{f}^{(t_{init})} = \left( f(\mathbf{x}^1, \mathbf{u}^1), \ldots , f(\mathbf{x}^{t_{init}}, \mathbf{u}^{t_{init}}) \right)$
\State $\mathbf{g}_1^{(t_{init})} = \left( g_1(\mathbf{x}^1, \mathbf{u}^1), \ldots , g_1(\mathbf{x}^{t_{init}}, \mathbf{u}^{t_{init}}) \right)$ 
\State \vdots
\State $\mathbf{g}_l^{(t_{init})} = \left( g_l (\mathbf{x}^1, \mathbf{u}^1), \ldots , g_l(\mathbf{x}^{t_{init}}, \mathbf{u}^{t_{init}}) \right)$ 
\State $t \gets t_{init}$
\While{$t \leq t_{max}$}
\State Train the GPs $F$ and $G_p$ for $p = 1, \ldots , l$ w.r.t. $(\mathcal{D}_f^{(t)},\mathbf{f}^{(t)})$ and $(\mathcal{D}_{g_p}^{(t)},\mathbf{g}_p^{(t)})$
\State Compute the targeted variable by maximizing the EFI (Equation \eqref{eq:EFI3}),\\
\hspace{1cm} $\mathbf{x}_{targ} = \mbox{arg} \min_{\mathbf{x}} \mathbb E[FI^{(t)}(\mathbf{x})]$
\State Set $\mathbf{x^{t+1}} = \mathbf{x}_{targ}$
\State Compute the random parameters of the objective function by minimizing\\
\hspace{1cm} the one-step-ahead improvement variance (Equation \eqref{eq:S1}) \\
\hspace{1cm} $\mathbf{u}_f^{t+1} = \mbox{arg} \min_{\mathbf{u}} S_f(\mathbf x_{targ},\mathbf u)$
\State Determine the next constraint to evaluate and its random parameter \\
\hspace{1cm} by minimizing the one-step-ahead variance \\
\hspace{1cm} of the feasibility probability (Equation \eqref{SC2_ConstSel})\\
\hspace{1cm} $\{\mathbf u_g^{t+1}, p \} = \mbox{arg} \min_{\mathbf{u},p'} S_g(\mathbf x_{targ},\mathbf u,p')$
\State $\mathcal{D}_f^{(t)} \gets \{\mathcal{D}_f^{(t)}, (\mathbf{x}^{t+1}, \mathbf{u}^{t+1}_f)\} $ \quad , \quad
$\mathbf{f}^{(t)} \gets \{\mathbf{f}^{(t)}, f(\mathbf{x}^{t+1}, \mathbf{u}^{t+1}_f)\} $
\State $\mathcal{D}_{g_p}^{(t)} \gets \{\mathcal{D}_{g_p}^{(t)}, (\mathbf{x}^{t+1}, \mathbf{u}^{t+1}_g)\} $ \quad , \quad
$\mathbf{g}_p^{(t)} \gets \{\mathbf{g}_p^{(t)}, g_p(\mathbf{x}^{t+1}, \mathbf{u}^{t+1}_g)\} $
\State $t \gets t+1$ \quad , \quad $t_p \gets t_p+1$
\EndWhile \\
\Return the best feasible point according to the means of the processes,\\
\hspace{1cm} $\arg \min_{\mathbf{x}} m^{(t)}_Z(\mathbf{x}) ~\mbox{s.t.}~ \mathbb E[ C^{(t)}(\mathbf{x}) ] \leq 0$ 
\end{algorithmic}
\end{algorithm}

\begin{algorithm}
\caption{MMCU algorithm}\label{alg:MMCU}
\begin{algorithmic}
\State Initialize $t_{init}$, $t_{max}$
\State Define the initial data set $\mathcal{D}^{(t_{init})} = \{ (\mathbf{x}^i,\mathbf{u}^i), \ i= 1,\ldots t_{init} \}$
\State Evaluate the objective and the constraint functions : 
\State $\mathbf{f}^{(t_{init})} = \left( f(\mathbf{x}^1, \mathbf{u}^1), \ldots , f(\mathbf{x}^{t_{init}}, \mathbf{u}^{t_{init}}) \right)$
\State $\mathbf{g}_1^{(t_{init})} = \left( g_1(\mathbf{x}^1, \mathbf{u}^1), \ldots , g_1(\mathbf{x}^{t_{init}}, \mathbf{u}^{t_{init}}) \right)$ 
\State \vdots
\State $\mathbf{g}_l^{(t_{init})} = \left( g_l (\mathbf{x}^1, \mathbf{u}^1), \ldots , g_l(\mathbf{x}^{t_{init}}, \mathbf{u}^{t_{init}}) \right)$ 
\State $t \gets t_{init}$
\While{$t \leq t_{max}$}
\State Train the GPs $F$ w.r.t. $\left(\mathcal{D}^{(t)}, \mathbf{f}^{(t)}\right)$ and \\
\hspace{1cm} and $\mathbf{G}$ w.r.t. $\left(\underbrace{\{ \mathcal{D}^{(t)},..., \mathcal{D}^{(t)}\}}_{l\text{ times}}, \{\mathbf{g}_1^{(t)}, \ldots, \mathbf{g}_l^{(t)}\}\right)$
\State Compute the targeted variable by maximizing the EFI (Equation \eqref{eq:EFI3}),\\
\hspace{1cm} $\mathbf{x}_{targ} = \mbox{arg} \min_{\mathbf{x}} \mathbb E[FI^{(t)}(\mathbf{x})]$
\State Set $\mathbf{x^{t+1}} = \mathbf{x}_{targ}$
\State Compute the random parameter $\mathbf u^{t+1}$ as the minimizer \\
\hspace{1cm} of the proxy to the one-step-ahead EFI variance (Equation \eqref{eq:proxy}), \\
\hspace{1cm} $\mathbf{u}^{t+1} = \arg \min_\mathbf{u} S(\mathbf{x}_{targ},\mathbf{u}) $
\State $\mathcal{D}^{(t)} \gets \{\mathcal{D}^{(t)}, (\mathbf{x}^{t+1}, \mathbf{u}^{t+1})\} $
\State $\mathbf{f}^{(t)} \gets \{\mathbf{f}^{(t)}, f(\mathbf{x}^{t+1}, \mathbf{u}^{t+1})\} $
\State $\mathbf{g}_p^{(t)} \gets \{\mathbf{g}_p^{(t)}, g_p(\mathbf{x}^{t+1}, \mathbf{u}^{t+1})\} $ for $p = 1, \ldots, l$
\State $t \gets t+1$
\EndWhile\\
\Return the best feasible point according to the means of the processes,\\
\hspace{1cm} $\arg \min_{\mathbf{x}} m^{(t)}_Z(\mathbf{x}) ~\mbox{s.t.}~ \mathbb E[ C^{(t)}(\mathbf{x}) ] \leq 0$ 
\end{algorithmic}
\end{algorithm}

\begin{algorithm}
\caption{MMCS algorithm}\label{alg:MMCS}
\begin{algorithmic}
\State Initialize $t_{init}$, $t_{max}$
\State Define the initial data set $\mathcal{D}_f^{(t_{init})} = \{ (\mathbf{x}^i,\mathbf{u}^i), \ i = 1,\ldots t _{init}\}$ \\
\hspace{1cm} and $\mathcal{D}_{\mathbf g}^{(t_{init})} = \{ (\mathbf{x}^i,\mathbf{u}^i,p) ~,~i = 1,\ldots t_{init}~,~p=1,\ldots,l \}$
\State Evaluate the objective and the constraint functions : 
\State $\mathbf{f}^{(t_{init})} = \left( f(\mathbf{x}^1, \mathbf{u}^1), \ldots , f(\mathbf{x}^{t_{init}}, \mathbf{u}^{t_{init}}) \right)$
\State $\mathbf g^{(t_{init})} = \left( g_1(\mathbf{x}^1, \mathbf{u}^1), \ldots , g_1(\mathbf{x}^{t_{init}}, \mathbf{u}^{t_{init}}), \ldots ,
g_l (\mathbf{x}^1, \mathbf{u}^1), \ldots , g_l(\mathbf{x}^{t_{init}}, \mathbf{u}^{t_{init}}) \right)$
\State $t \gets t_{init}$
\While{$t \leq t_{max}$}
\State Train the GPs $F$ w.r.t. $(\mathcal{D}_f^{(t)},\mathbf{f}^{(t)})$  and $\mathbf{G}$ w.r.t. $(\mathcal{D}_{\mathbf g}^{(t)} , \mathbf{g}^{(t)})$
\State Compute the targeted variable by maximizing the EFI (Equation \eqref{eq:EFI3}),\\
\hspace{1cm} $\mathbf{x}_{targ} = \mbox{arg} \min_{\mathbf{x}} \mathbb E[FI^{(t)}(\mathbf{x})]$
\State Set $\mathbf{x^{t+1}} = \mathbf{x}_{targ}$
\State Compute the random parameters of the objective function by minimizing\\
\hspace{1cm} the one-step-ahead improvement variance (Equation \eqref{eq:S1}) \\
\hspace{1cm} $\mathbf{u}_f^{t+1} = \mbox{arg} \min_{\mathbf{u}} S_f(\mathbf x_{targ},\mathbf u)$
\State Determine the next constraint to evaluate and its random parameter \\
\hspace{1cm} by minimizing the one-step-ahead variance \\
\hspace{1cm} of the feasibility probability (Equation \eqref{SC2_ConstSel})\\
\hspace{1cm} $\{\mathbf u_g^{t+1}, p \} = \mbox{arg} \min_{\mathbf{u},p'} S_g(\mathbf x_{targ},\mathbf u,p')$

\State $\mathcal{D}_f^{(t)} \gets \{\mathcal{D}_f^{(t)}, (\mathbf{x}^{t+1}, \mathbf{u}^{t+1}_f)\} $ ~,~
$\mathbf{f}^{(t)} \gets \{\mathbf{f}^{(t)}, f(\mathbf{x}^{t+1}, \mathbf{u}^{t+1}_f)\} $
\State $\mathcal{D}_{\mathbf g}^{(t)} \gets \{\mathcal{D}_{\mathbf g}^{(t)}, (\mathbf{x}^{t+1}, \mathbf{u}_g^{t+1},p)\} $ ~,~
$\mathbf{g}^{(t)} \gets \{\mathbf{g}^{(t)}, g_p(\mathbf{x}^{t+1}, \mathbf{u}^{t+1}_g)\} $ 
\State $t \gets t+1$
\EndWhile\\
\Return the best feasible point according to the means of the processes,\\
\hspace{1cm} $\arg \min_{\mathbf{x}} m^{(t)}_Z(\mathbf{x}) ~\mbox{s.t.}~ \mathbb E[ C^{(t)}(\mathbf{x}) ] \leq 0$ 
\end{algorithmic}
\end{algorithm}

\end{document}